\documentclass[a4paper,12pt]{article}
\usepackage[english]{babel}   
\usepackage{indentfirst}
\usepackage[T2A]{fontenc}
\usepackage[utf8]{inputenc}

\renewcommand{\epsilon}{\ensuremath{\varepsilon}}
\renewcommand{\phi}{\ensuremath{\varphi}}
\renewcommand{\kappa}{\ensuremath{\varkappa}}
\renewcommand{\le}{\ensuremath{\leqslant}}

\renewcommand{\ge}{\ensuremath{\geqslant}}

\usepackage{amsmath,amsfonts,amssymb,amsthm,mathtools} 
\usepackage{icomma} 

\usepackage{graphicx}  
\graphicspath{{images/}{images2/}}  
\usepackage{geometry} 
\usepackage{setspace} 
\usepackage{lastpage} 
\usepackage{soul} 
  \usepackage{hyperref}
\usepackage[usenames,dvipsnames,svgnames,table,rgb]{xcolor}
\begin{document}

\noindent  \begin{center} {\huge About one matrix of composite numbers\begin{spacing}{1.7}\end{spacing} and her applications} \end{center}

\textbf{}

\begin{center} {\large Garipov Ilshat Ilsurovich }\end{center}
\begin{spacing}{-1}\end{spacing}
\begin{center} \textit {\normalsize Russia, Republic of Tatarstan, Naberezhnye Chelny}\end{center}
\begin{spacing}{-1}\end{spacing}
\begin{center} \textit {\normalsize e-mail: mathsciencegaripovii@gmail.com.}\end{center}
\renewcommand{\abstractname}{Abstract}
\begin{abstract}
{\small Scientific paper is devoted to research of \textit{T}-matrix -- matrix of composite numbers $6h\pm 1$ in special view, and to her application in number theory.}


Results of this paper were presented in XVIII International Conference «Algebra, Number Theory and Discrete Geometry: modern problems, applications and problems of history».
\end{abstract}

\textbf{List of symbols}

${\rm \mathbb{N} } $ -- set of all natural numbers.

${\rm \mathbb{N} } _{0} $ -- set of all natural numbers with zero. 

${\rm \mathbb{Z} } $ -- set of all integers.

${\rm   \mathbb{P} } $ -- set of all prime numbers.

${\rm  \mathbb{R} }$-- set of all real numbers.

$T$ -- matrix comprising all defining and not defining elements.

 $T^*$ -- matrix obtained from $T$-matrix by deleting first row.

$\widetilde{T}$ -- set of all elements of $T$-matrix.

$D_{T} $ -- set of all defining elements of $T$-matrix.

$nD_{T} $ -- set of all not defining elements of $T$-matrix.

${\rm M} _{T} $ -- set of all leading elements of $T$-matrix.

$D_{T_{k} } $ -- set of all defining elements from $k$-row of $T$-matrix, $k\in {\rm  \mathbb{N}} $.

$nD_{T_{k} } $-- set of all not defining elements from $k$-row of $T$-matrix, $k\in {\rm  \mathbb{N}} $.

$D_{T,6h+1} $ -- set of all $T$-matrix defining elements of the form $6h+1$.

$D_{T,6h-1} $ -- set of all $T$-matrix defining elements of the form $6h-1$.

$nD_{T,6h+1} $ -- set of all $T$-matrix not defining elements of the form $6h+1$.

$nD_{T,6h-1} $ -- set of all $T$-matrix not defining elements of the form $6h-1$.

$\pi (x)$ -- function counting the number of prime numbers less than or equal to $x\in {\rm \mathbb{R}}$.

$\pi _{{\rm M} _{T} } (x)$ -- function counting the number of $T$-matrix leading elements less than or equal to $ x\in {\rm \mathbb{R}}$.

\begin{center}
\textbf{\large Introduction. 1. Construct \textit{T}-matrix}
\end{center}
\begin{spacing}{0.5}\end{spacing}

We first construct numerical sequences $T_{k} \equiv \left(a(k;n)\right)_{n=1}^{\infty } $, where $k$ is a number of sequence $T_{k}$, $k\in {\rm \mathbb{N}}$; $a(k;n)$ is the $n$-th element of sequence $T_{k}$.

The elements $a(k;n)$ are defined as follows:
\begin{spacing}{0.8}\end{spacing}
\begin{equation} \label{(1)}
a(k;n)\equiv p(k)\cdot \left(5+2\cdot \left\lfloor \frac{n}{2} \right\rfloor +4\cdot \left\lfloor \frac{n-1}{2} \right\rfloor \right),          
\end{equation} 
where $p(k)$ is the $k$-th element of sequence $\left(p(k)\right)_{k=1}^{\infty } $ of prime numbers. The elements $p(k)$ are defined under rule:
\begin{equation}  \label{(2)}
p(k)\equiv p_{k+2} ,           
\end{equation} 
\begin{spacing}{1.1}\end{spacing}
\noindent where $p_{i}$ is the $i$-th prime number in sequence of all prime numbers.                   

As a result, we have a matrix $T\equiv \left(a(k;n)\right)_{\infty \times \infty } $ in which $k$-row presents all elements of sequence $T_{k} $. In this case, $a(k;n)$ is an element of matrix $T$, located in $k$-th row and $n$-th column (see Table №1).

\begin{center} \textbf{Table №1. }Fragment submatrix in matrix $T$ for the case of $k=\overline{1;15},{\rm \; }n=\overline{1;30}$\end{center}
 \begin{spacing}{-1.3}\end{spacing}\begin{center}\includegraphics*[scale=0.55]{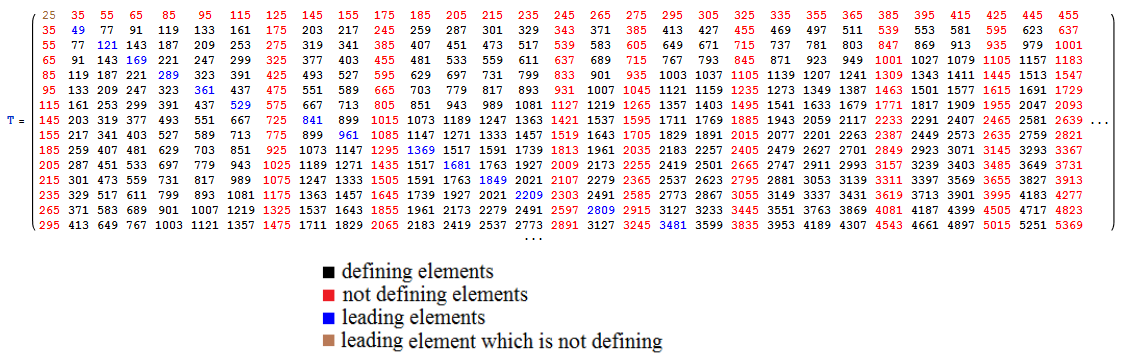} \end{center}

THEOREM 1.1. \begin{equation} \label{(3)} \left(\forall k,n\in {\rm \mathbb{N}} \right)\left(a(k;n)=p(k)\cdot \left(3n+\frac{3-(-1)^{n} }{2} \right)\right).\end{equation} 

PROOF. Applying Hermite's identity (see [1]):
\begin{equation} \label{(4)} 
\left(\forall m\in {\rm \mathbb{N}} \right)\left(\left\lfloor mx\right\rfloor =\sum _{i=0}^{m-1}\left\lfloor x+\frac{i}{m} \right\rfloor  \right).  
\end{equation} 

If $m=2$, then identity \eqref{(4)}  has the form: \begin{spacing}{0.8}\end{spacing}\begin{equation} \label{(5)} \left\lfloor 2x\right\rfloor =\left\lfloor x\right\rfloor +\left\lfloor x+\frac{1}{2} \right\rfloor .\end{equation} \begin{spacing}{1.1}\end{spacing}
 
For all $n\in {\rm \mathbb{N}} $ we substitute in \eqref{(5)} $x=\frac{n-1}{2} $. Then,
\[{\rm \; \; \; \; \; }\left\lfloor 2\cdot \left(\frac{n-1}{2} \right)\right\rfloor =\left\lfloor \frac{n-1}{2} \right\rfloor +\left\lfloor \frac{n-1}{2} +\frac{1}{2} \right\rfloor {\rm \;}\Leftrightarrow {\rm \;}\left\lfloor n-1\right\rfloor =\left\lfloor \frac{n-1}{2} \right\rfloor +\left\lfloor \frac{n-1+1}{2}  \right\rfloor {\rm \;}\mathop{\Rightarrow }\limits^{n-1\in {\rm \mathbb{Z}} } {\rm \; }\] 
\begin{equation} \label{(6)} 
 \mathop{\Rightarrow } {\rm \; \; } n-1=\left\lfloor \frac{n-1}{2} \right\rfloor +\left\lfloor \frac{n}{2} \right\rfloor {\rm \;}\Leftrightarrow {\rm \;}\left\lfloor \frac{n-1}{2} \right\rfloor =n-1-\left\lfloor \frac{n}{2} \right\rfloor . 
\end{equation} 

If to compare the statement of Theorem 1.1 and \eqref{(1)}, then it would suffice to prove:
\[5+2\cdot \left\lfloor \frac{n}{2} \right\rfloor +4\cdot \left\lfloor \frac{n-1}{2} \right\rfloor =3n+\frac{3-(-1)^{n} }{2} .\] 
\[5+2\cdot \left\lfloor \frac{n}{2} \right\rfloor +4\cdot \left\lfloor \frac{n-1}{2} \right\rfloor \mathop{=}\limits^{\eqref{(6)} } 5+2\cdot \left\lfloor \frac{n}{2} \right\rfloor +4\cdot \left(n-1-\left\lfloor \frac{n}{2} \right\rfloor \right)=\]
\[=5+2\cdot \left\lfloor \frac{n}{2} \right\rfloor +4n-4-4\cdot \left\lfloor \frac{n}{2} \right\rfloor =4n+1-2\cdot \left\lfloor \frac{n}{2} \right\rfloor .\] 
\[\left\lfloor \frac{n}{2} \right\rfloor =\frac{n}{2} -\left\{\frac{n}{2} \right\}=\frac{n}{2} -\frac{1-(-1)^{n} }{4} \Rightarrow \] 
\[\Rightarrow {\rm \; \; }4n+1-2\cdot \left\lfloor \frac{n}{2} \right\rfloor =4n+1-2\cdot \left(\frac{n}{2} -\frac{1-(-1)^{n} }{4} \right)=\]
\[=4n+1-n+\frac{1-(-1)^{n} }{2} =3n+\frac{3-(-1)^{n} }{2} .\] 

Theorem 1.1 is proved.

DEFINITION 1.1. An element $a(k;n)$ of matrix $T$ is called defining if:

1)  $a(k;n)$ is not divisible by 5;

2) $a(k;n)$ can be expressed as a product of some two prime numbers, that is:
\begin{spacing}{0.5}\end{spacing}\begin{equation} \label{(7)} 
5\not |a(k;n){\rm \; \; }\wedge {\rm \; \; }\left(\exists p_{1} ,p_{2} \in {\rm \mathbb{P}} \right)\left(a(k;n)=p_{1} \cdot p_{2} \right).  
\end{equation} \begin{spacing}{1.1}\end{spacing}

DEFINITION 1.2. An element $a(k;n)$ of matrix $T$ is called not defining if he does not satisfy condition \eqref{(7)}. 

DEFINITION 1.3. An element $a(k;n)$ of matrix $T$ is called leading if:
\begin{spacing}{0.5}\end{spacing}\begin{equation} \label{(8)} 
a(k;n)=p^{2} (k).    
\end{equation} \begin{spacing}{1.1}\end{spacing}

DEFINITION 1.4. $T$-matrix is called matrix comprising all defining and not defining elements.

COMMENT. One of the factors in factorization of any element $a(k;n)$ of $T$-matrix always the same as the number $p(k)$ defined by rule \eqref{(2)}.  
\begin{spacing}{1.7}\end{spacing}
LEMMA 1.1. $\left(3n+\frac{3-(-1)^{n} }{2} \right)_{n=1}^{\infty } $ is a sequence of all numbers of the form $6h\pm 1$:
\begin{spacing}{0.5}\end{spacing}\begin{equation} \label{(9)} 
5;{\rm \; }7;{\rm \; }11;{\rm \; }13;{\rm \; }17;{\rm \; }19;{\rm \; }23;{\rm \; }25;{\rm \; }...{\rm \; };{\rm \; }6h-1;{\rm \; 6}h+{\rm 1;...}.  
\end{equation} \begin{spacing}{1.1}\end{spacing}

PROOF. $n=2h-1, h\in {\rm \mathbb{N}} {\rm \; \; }\Leftrightarrow {\rm \; \; }3n+\frac{3-(-1)^{n} }{2} ={\rm \; }3\cdot (2h-1)+\frac{3-(-1)^{2h-1} }{2}=6h-3+2=6h-1$. \[ n=2h, h\in {\rm \mathbb{N}} {\rm \; \; }\Leftrightarrow {\rm \; \; }3n+\frac{3-(-1)^{n} }{2} ={\rm \; }3\cdot (2h)+\frac{3-(-1)^{2h} }{2} =6h+1.\] 

 As a result, we will have sequence of form \eqref{(9)}. Lemma 1.1 is proved.
\begin{spacing}{1.7}\end{spacing}
Then let $f(n)\equiv 3n+\frac{3-(-1)^{n} }{2} $.

\begin{center}
\textbf{\large 2. The simplest properties of \textit{T}-matrix }
\end{center}

\begin{spacing}{0.5}\end{spacing}
PROPERTY 2.1. $T$-matrix is matrix of infinite order.

PROOF. It follows from the definition of sequence $T_{k} $ for all $k\in {\rm \mathbb{N}} $, and  \eqref{(1)} that the $T$-matrix consists of an infinite number of rows and columns. It would mean that $T$-matrix has infinite order. Property 2.1 is proved.  

PROPERTY 2.2. $rang(T)=1$.

PROOF. Using presentation of $T$-matrix elements in Theorem 1.1, we find the highest of orders of every possible nonzero minors of $T$-matrix. This order will determine rank of $T$- matrix. Find rank $rang(T)$ of $T$-matrix with the method of bordering minors (about the method see [3]): 

$1) \left|a(1;1)\right|=\left|25\right|=25\ne 0.$

2) For any $k,n\in {\rm \mathbb{N}} \backslash \{ 1\} $ we have:
\begin{spacing}{-0.2}\end{spacing}\[\left|\begin{array}{cc} {a(1;1)} & {a(1;n)} \\ {a(k;1)} & {a(k;n)} \end{array}\right|{\rm \; }\mathop{=}\limits^{\eqref{(3)}} {\rm \; }\left|\begin{array}{cc} {p(1)\cdot f(1)} & {p(1)\cdot f(n)} \\ {p(k)\cdot f(1)} & {p(k)\cdot f(n)} \end{array}\right|=f(1)\cdot \left|\begin{array}{cc} {p(1)} & {p(1)\cdot f(n)} \\ {p(k)} & {p(k)\cdot f(n)} \end{array}\right|=\]
\[=f(1)\cdot f(n)\cdot \left|\begin{array}{cc} {p(1)} & {p(1)} \\ {p(k)} & {p(k)} \end{array}\right|=f(1)\cdot f(n)\cdot {\rm \; }0={\rm 0}.\] 

Therefore, $rang(T)=1$. Property 2.2 is proved.   

THEOREM 2.1. If an infinite matrix $A$ has nonzero rank, then there exists matrix $A$-induced infinite determinant, and his value is zero (see [2]).

PROPERTY 2.3. $\det (T)=0$.\textbf{}

PROOF. By Property 2.1, $T$-matrix is an infinite. It follows from Property 2.2 that $T$-matrix has nonzero rank. Using Theorem 2.1, we get that matrix $T$-induced infinite determinant $\det (T)$ equals 0. Property 2.3 is proved.

PROPERTY 2.4. $\left(\forall n\in {\rm \mathbb{N}} \right)\left(a(1;n)\in nD_{T} \right){\rm \; \; }\wedge {\rm \; \; }a(1;1)\in {\rm M} _{T} $.

PROOF. By construction: $T_{1} \equiv \left(a(1;n)\right)_{n=1}^{\infty } $.  Let us use \eqref{(3)} for presentation of elements $a(1;n)$ of sequence $T_{1} $ and Definition 1.2. 
\begin{spacing}{0.5}\end{spacing} 
\[(\forall n\in {\rm \mathbb{N}} )(a(1;n)=p(1)\cdot f(n)=5\cdot f(n)){\rm  \; }\Rightarrow {\rm  \; }(\forall n\in {\rm \mathbb{N}} )(5|{\rm \; }a(1;n)){\rm  \; }\Rightarrow\]
\[\Rightarrow {\rm  \; }(\forall n\in {\rm \mathbb{N}})(a(1;n)\in nD_{T_{1} }). \]\begin{spacing}{0.8}\end{spacing}

So, all elements of sequence $T_{1}$ are not defining. 

Let us show that element $a(1;1)$ is leading. Using Definition 1.3, we get:
\begin{center}$a(1;1){\rm \; }\mathop{=}\limits^{\eqref{(3)}} {\rm \; }p(1)\cdot f(1){\rm \; \; }\mathop{=}\limits^{p(1)=f(1)} {\rm \; }p^{2} (1){\rm \; \; }\Rightarrow {\rm \; \; }a(1,1)\in {\rm M} _{T} .$\end{center} \begin{spacing}{0.8}\end{spacing}

Property 2.4 is proved.

PROPERTY 2.5. There exists a leading element of set ${\rm M} _{T} $ that do not lie on the main diagonal of $T$-matrix.

PROOF. Suppose otherwise: all leading elements of set ${\rm M} _{T} $ lie on the main diagonal of $T$-matrix. Then with respect to the assumption and \eqref{(3)}:
\begin{spacing}{0.5}\end{spacing} \[a(k;k)=p^{2} (k){\rm \; \; }\Leftrightarrow {\rm \; \; }p(k)\cdot f(k)=p^{2} (k){\rm \; \; }\Leftrightarrow {\rm \; \; }p(k)=f(k).\] \begin{spacing}{1.1}\end{spacing}

Since $k$ is any natural number, choose a number $k=8$:
\begin{spacing}{-0.2}\end{spacing}\[p(8)\mathop{=}\limits^{\eqref{(2)}} p_{10} {\rm \; \; }\Leftrightarrow {\rm \; \; }f(8)=p_{10} {\rm \; \; }\Leftrightarrow {\rm \; \; }3\cdot 8+\frac{3-(-1)^{8} }{2} =29{\rm \; \; }\Leftrightarrow {\rm \; \; 25}=29{\rm \; \; }\Leftrightarrow {\rm \; \; } \varnothing .\]\begin{spacing}{1.1}\end{spacing}

As a result, a contradiction. Property 2.5 is proved.

PROPERTY 2.6. $D_{T} \cap nD_{T} =\varnothing $.\textbf{}

PROOF. Given Definition 1.1 and Definition 1.2 for defining and not defining elements of $T$-matrix, we get:
\begin{spacing}{0.5}\end{spacing}
\begin{equation} \label{(10)} 
D_{T} =\left\{a\in \widetilde{T}{\rm \; |\; }\left(\exists p_{1} ,p_{2} \in {\rm  \mathbb{P}} \right) \left(a=p_{1} \cdot p_{2} \right){\rm \; \; }\wedge {\rm \; \; }5\not |a\right\};  
\end{equation} 
\begin{equation} \label{(11)} 
nD_{T} =\left\{a\in \widetilde{T}{\rm \; |\; }\neg \left(\left(\exists p_{1} ,p_{2} \in {\rm \mathbb{P}} \right)\left(a=p_{1} \cdot p_{2} \right){\rm \; \; }\wedge {\rm \; \; }5\not |a\right)\right\}.  
\end{equation} 
\[D_{T} \cap nD_{T} \mathop{=}\limits^{\eqref{(10)},{\rm \; \eqref{(11)}}} {\rm \; \; }\left\{a\in \widetilde{T}{\rm \; |\; }\left(\exists p_{1} ,p_{2} \in {\rm \mathbb{P}} \right) \left(a=p_{1} \cdot p_{2} \right){\rm \; \; }\wedge {\rm \; \; }5\not |a \right\}{\rm \; }\cap \]
 \[\cap {\rm \; }\left\{a\in \widetilde{T}{\; | \;}\neg \left(\left(\exists p_{1} ,p_{2} \in {\rm\mathbb{P}} \right)\left(a=p_{1} \cdot p_{2} \right){\rm \; \; }\wedge {\rm \; \; }5\not |a\right)\right\}=\] 
\[=\left\{a\in \widetilde{T}{|} \left(\exists p_{1} ,p_{2} \in {\rm \mathbb{P}} \right)\left(a=p_{1} \cdot p_{2} \right){\rm \;}\wedge {\rm \; }5\not |a{\rm \;}\wedge {\rm \;}\neg \left(\left(\exists p_{1} ,p_{2} \in {\rm \mathbb{P}} \right)\left(a=p_{1} \cdot p_{2} \right){\rm \;}\wedge {\rm \;}5\not |a\right)\right\}= \varnothing .\] 

Property 2.6 is proved.\textbf{}

PROPERTY 2.7. $M_{T} \cap nD_{T} =\{ 25\} $.

PROOF. Given Definition 1.3 for leading element of $T$-matrix
\begin{spacing}{0.5}\end{spacing}\begin{equation} \label{(12)} 
{\rm M} _{T} =\left\{a\in \widetilde{T}{\rm \; |\; }\left(\exists p\in {\rm \mathbb{P}} \right)\left(a=p^{2} \right)\right\}.  
\end{equation} \begin{spacing}{1.1}\end{spacing}
It follows from \eqref{(11)} and \eqref{(12)} that 
\[{\rm M} _{T} \cap nD_{T}=\left\{a\in \widetilde{T}{\rm \; |\; }\left(\exists p\in {\rm  \mathbb{P}} \right)\left(a=p^{2} \right)\right\} \cap \]
\begin{spacing}{-0.3}\end{spacing}
\[ \cap \left\{a\in \widetilde{T}{\rm \; |\; \; }\neg \left(\left(\exists p_{1} ,p_{2} \in {\rm  \mathbb{P}} \right)\left(a=p_{1} \cdot p_{2} \right){\rm \;}\wedge {\rm \;}5\not |a\right)\right\}=\] 
\begin{spacing}{0.3}\end{spacing}
\[=\left\{a\in \widetilde{T}{\rm \; |\; }\left(\exists p\in {\rm \mathbb{P}} \right)\left(a=p^{2} \right){\rm \;}\wedge {\rm \;}\neg \left(\left(\exists p_{1} ,p_{2} \in {\rm  \mathbb{P}} \right)\left(a=p_{1} \cdot p_{2} \right){\rm \; \; }\wedge {\rm \;  \; }5\not |a\right)\right\}\equiv \] 
\[\equiv \left\{a\in \widetilde{T}{\rm \; |\; }\left(\exists p\in {\rm  \mathbb{P}} \right)\left(a=p^{2} \right){\rm \; }\wedge {\rm \; }\left(\neg \left(\exists p_{1} ,p_{2} \in {\rm \mathbb{P}} \right) \left(a=p_{1} \cdot p_{2} \right){\rm \; \; }\vee {\rm \; \; }\neg \left(5\not |a\right)\right) \right\}\equiv \] 
\[\equiv \left\{{\rm \; }a\in \widetilde{T}{\rm \; |\; }\left(\exists p\in {\rm  \mathbb{P}} \right)\left(a=p^{2} \right){\rm \; }\wedge {\rm \;}\left(\left(\forall p_{1} ,p_{2} \in {\rm \mathbb{P}} \right)\left(a\ne p_{1} \cdot p_{2} \right){\rm \; \; }\vee {\rm \; \; }5|a\right) \right\}\equiv \] 
\[\equiv \left\{a\in \widetilde{T}{\rm \; |\; }\left(\exists p\in {\rm \mathbb{P}} \right) \left(a=p^{2} \right){\rm \; }\wedge {\rm \; }\left(\forall p_{1} ,p_{2} \in {\rm  \mathbb{P}} \right) \left(a\ne p_{1} \cdot p_{2} \right){\rm \; }\vee {\rm \; }\left(\exists p\in {\rm  \mathbb{P}} \right)\left(a=p^{2} \right){\rm \; }\wedge {\rm \; }5|a \right\} \equiv \] 
\[\equiv \left\{a\in \widetilde{T}{\rm \; |\; }\left(\exists p\in {\rm \mathbb{P}} \right) \left(a=p^{2} \right){\rm \;}\wedge {\rm \;}5|a \right\}=\{ 25\} .\] 

Property 2.7 is proved.

PROPERTY 2.8.The sequence $\left(p^{2} (k)\right)_{k=1}^{\infty } $ of $T$-matrix leading elements is ascending.

This property follows from the fact that the sequence $\left(p(k)\right)_{k=1}^{\infty } $ is ascending.

DEFINITION 2.1. A semiprime number is a composite number that is the product of two (possibly equal) primes.

PROPERTY 2.9. All defining elements of $T$-matrix are semiprimes.

PROOF. Choose any defining element $a(k;n)$ of $T$-matrix. Then it follows from \eqref{(7)} that
\begin{spacing}{0.5}\end{spacing}
\begin{equation} \label{(13)} 
\left(\exists p_{1} ,p_{2} \in {\rm \mathbb{P}} \right)\left(a(k;n)=p_{1} \cdot p_{2} \right).  
\end{equation} 
In this case, given Definition 2.1, \eqref{(13)} will mean that $a(k;n)$ is semiprime number. 

Property 2.9 is proved.

COMMENT. Not every semiprime number is defining element of $T$-matrix. For example,

\noindent semiprime number 25 isn't defining element of $T$-matrix, since condition 1) of Definition 1.1 is violated. 

\begin{center}
\textbf{\large 3. Support properties and theorems from number theory and set theory}
\end{center}

\begin{spacing}{0.9}\end{spacing}
PROPERTY 3.1 (of prime numbers). 
\begin{spacing}{0.5}\end{spacing}
\[\left(\forall p\in {\rm \mathbb{P}} \right)\left(p>3 {\rm \;}\Rightarrow \left(\exists t\in {\rm\mathbb{N}} \right) \left(p=6t+1{\rm \; }\vee {\rm \; }p=6t-1\right)\right).\]

THEOREM 3.2. The set of prime numbers of the form $6t+1$ is infinite.

THEOREM 3.3. The set of prime numbers of the form $6t-1$ is infinite.

THEOREM 3.4 (Euclid's theorem). The set of prime numbers is infinite (see [4]).

Also we need some theorems from set theory.

THEOREM 3.5. If to remove finite subset from countable set, then the remaining set is countable.

THEOREM 3.6. The union of countably many countable sets is countable.

THEOREM  3.7. The union of a finite number of countable sets is countable.

THEOREM 3.8. A subset of a countable set is not more than countable.

THEOREM 3.9. The union of not more than countable set and countable set is countable (see [5]).

Further, we will formulate the fundamental theorem of arithmetic for natural numbers.  

THEOREM 3.10. Every positive integer except the number 1 can be represented in exactly one way apart from rearrangement as a product of one or more primes (see [4]).

\begin{center}
\textbf{\large 4. The basic theorems about elements of \textit{T}-matrix}
\end{center}

\begin{spacing}{0.5}\end{spacing}
LEMMA 4.1.
\begin{spacing}{0.1}\end{spacing}
\[\left(\forall b\in \widetilde{T}\right)\left(\exists h\in {\rm \mathbb{N}} \right)\left(b=6h+1{\rm \; \; }\vee {\rm \; \; }b=6h-1\right).\]
\begin{spacing}{1.2}\end{spacing}
PROOF. Choose any element $a(k;n)$ of $T$-matrix. Using Property 3.1, let us consider 2 cases.

 \textbf{Case 1.} For some number $k\in {\rm \mathbb{N}} $:
 \begin{spacing}{0.2}\end{spacing}
 \begin{equation} \label{(14)} 
 \left(\exists t\in {\rm \mathbb{N}} \right)\left(p(k)=6t+1\right).
  \end{equation}
  
Once we find this number $t$, given Lemma 1.1,
\begin{spacing}{0.8}\end{spacing}
\[a(k;n)\mathop{=}\limits^{\eqref{(3)}} p(k)\cdot f(n){\rm \; }\mathop{=}\limits^{\eqref{(14)}} \left(6t+1\right)\cdot f(n)=\left(6t+1\right)\cdot \left(6h_{1} \pm 1\right)=\]
\[=6h_{1} \cdot \left(6t+1\right)\pm \left(6t+1\right)=6h_{1} \cdot \left(6t+1\right)\pm 6t\pm 1=6\cdot \left(h_{1} \cdot \left(6t+1\right)\pm t\right)\pm 1{\rm \; \; }\Rightarrow\]
\[\Rightarrow{\rm \; \; }\left(\exists h\in {\rm \mathbb{N}} \right)\left(a(k;n)=6h\pm 1\right). \] 

\textbf{Case 2. }For some number $k\in {\rm \mathbb{N}} $:
\begin{spacing}{0.3}\end{spacing}\begin{equation} \label{(15)} 
\left(\exists t\in {\rm\mathbb{N}} \right)\left(p(k)=6t-1\right).     
\end{equation} \begin{spacing}{0.8}\end{spacing}
Once we find this number $t$, given Lemma 1.1,
\begin{spacing}{0.8}\end{spacing}
\[a(k;n)\mathop{=}\limits^{\eqref{(3)}} p(k)\cdot f(n){\rm \; }\mathop{=}\limits^{\eqref{(15)}} \left(6t-1\right)\cdot f(n)=\left(6t-1\right)\cdot \left(6h_{2} \pm 1\right)=\]
 \[6h_{2} \cdot \left(6t-1\right)\pm \left(6t-1\right)=6h_{2} \cdot \left(6t-1\right)\pm 6t\mp 1=6\cdot \left(h_{2} \cdot \left(6t-1\right)\pm t\right)\mp 1{\rm \;}\Rightarrow\]
 \[\Rightarrow{\rm \;}\left(\exists h\in {\rm \mathbb{N}} \right)\left(a(k;n)=6h\mp 1\right).\]

Thus, when choosing any element of $T$-matrix he will have the form $6h\pm 1,{\rm \; }h\in {\rm \mathbb{N}}$. 

Lemma 4.1 is proved. 

Also note that paragraph 5 will provide a theorem on how all composite numbers of the form $6h\pm 1,{\rm \; }h\in {\rm \mathbb{N}} $ fully exhaust elements of $T$-matrix.

LEMMA 4.2. 
\[\left(\forall k\in {\rm \mathbb{N}} \right)\left(\exists !n\in {\rm \mathbb{N}} \right)\left(a(k;n)\in {\rm M} _{T} \right).\]

PROOF. \textbf{Existence.} Choose any $k$-row of $T$-matrix. Since by Property 3.1 all prime numbers have the form $6t\pm 1$, $t\in {\rm  \mathbb{N}}$; then, given Lemma 1.1, there exists a number $n\in {\rm \mathbb{N}} $ such that $p(k)=f(n)$.

 In this case,
\[a(k;n)\mathop{=}\limits^{\eqref{(3)}} p(k)\cdot f(n)=p^{2} (k){\rm \; \; \; }\mathop{\Rightarrow }\limits^{{\rm \eqref{(8)}}} {\rm \; \; \; }a(k;n)\in {\rm M} _{T} .\] 

 \textbf{Uniqueness. }Suppose:
\[\left(\forall k\in {\rm  \mathbb{N}} \right)\left(\exists n_{1} ,n_{2} \in {\rm \mathbb{N}} \right)\left(n_{1} \ne n_{2} {\rm \; \; }\wedge {\rm \; \; }a(k;n_{1} )=p^{2} (k){\rm \; \; }\wedge {\rm \; \; }a(k;n_{2} )=p^{2} (k)\right).         \] 
\[a(k;n_{1} )=p(k)\cdot f(n_{1} {\rm )\; \; }\wedge {\rm \; \; }a(k;n_{2} )=p(k)\cdot f(n_{2} {\rm )\; \; }\Leftrightarrow\]
\[ \Leftrightarrow {\rm \; \; }p(k)=f(n_{1} {\rm )\; \; }\wedge {\rm \; \; }p(k)=f(n_{2} {\rm )\;}\Rightarrow {\rm \; \; }f(n_{1} {\rm )\; }=f(n_{2} {\rm )\; \; \; }\Leftrightarrow {\rm \; \; \; }n_{1} =n_{2} .\] 

As a result, a contradiction. The uniqueness is established. Lemma 4.2 is proved.

THEOREM 4.3.
\begin{spacing}{0.2}\end{spacing}
\[\left(\forall k,n\in {\rm  \mathbb{N}} \right)\left(a(k;n)\in {\rm M} _{T} {\rm \; }\Leftrightarrow {\rm  \; }a(k;1)=a(1;n)\right).\]
\begin{spacing}{1.7}\end{spacing}
PROOF. $a(k;n)\in {\rm M} _{T} {\rm \; \; }\mathop{\Leftrightarrow }\limits^{{\rm \eqref{(8)}}} {\rm \; \; }a(k;n)=p^{2} (k){\rm \; \; \; }\mathop{\Leftrightarrow }\limits^{\eqref{(3)}} {\rm \; \; }p(k)\cdot f(n)=p^{2} (k){\rm \; \; }\Leftrightarrow $ 
\begin{spacing}{0.5}\end{spacing}\[\Leftrightarrow {\rm \; \; } p(k)=f(n) {\rm \; \; }\Leftrightarrow {\rm \; \; }5\cdot p(k)=5\cdot f(n){\rm \; \; }\Leftrightarrow {\rm \; \; }p(k)\cdot 5=5\cdot f(n){\rm \; \; }\mathop{\Leftrightarrow }\limits^{\eqref{(2)},{\rm \; }\eqref{(9)}}\]
\begin{spacing}{-0.2}\end{spacing}\[\Leftrightarrow{\rm \; \; }p(k)\cdot f(1)=p(1)\cdot f(n){\rm \; \; \; }\mathop{\Leftrightarrow }\limits^{\eqref{(3)}} {\rm \; \; \; }a(k;1)=a(1;n).\] 
\begin{spacing}{0.8}\end{spacing}

Theorem 4.3 is proved.

THEOREM 4.4 (about columns comprising not defining elements of $T$-matrix).
\begin{spacing}{-0.5}\end{spacing}\[\left(\forall n\in {\rm \mathbb{N}} \right)\left(\left(\exists k\in {\rm \mathbb{N}} \right)\left(k>1{\rm \; \; }\wedge {\rm \; \; }a(k;n)\in nD_{T} \right){\rm \; \;}\Leftrightarrow {\rm \; \; }\left(\forall r\in {\rm \mathbb{N}} \right)\left(a(r;n)\in nD_{T} \right)\right).\]
\begin{spacing}{0.5}\end{spacing}
PROOF. \textbf{Necessity. } Suppose that opposite is true: 
\begin{spacing}{-0.5}\end{spacing}\[\neg \left(\forall n\in {\rm \mathbb{N}}\right)\left(\left(\exists k\in {\rm \mathbb{N}} \right)\left(k>1{\rm \; \; }\wedge {\rm \; \; }a(k;n)\in nD_{T} \right){\rm \;}\Rightarrow {\rm \;}\left(\forall r\in {\rm \mathbb{N}} \right)\left(a(r;n)\in nD_{T} \right)\right){\rm \; }\Leftrightarrow \] 
\begin{spacing}{-1.5}\end{spacing}\[\Leftrightarrow {\rm \;}\left(\exists n\in {\rm \mathbb{N}} \right)\left(\neg \left(\left(\exists k\in {\rm \mathbb{N}} \right)\left(k>1{\rm \; \; }\wedge {\rm \; \; }a(k;n)\in nD_{T} \right){\rm \; }\Rightarrow {\rm \; } \left(\forall r\in {\rm \mathbb{N}} \right)\left(a(r;n)\in nD_{T} \right)\right)\right){\rm \;}\Leftrightarrow \] 
\begin{spacing}{-1.5}\end{spacing}\[\Leftrightarrow {\rm \;}\left(\exists n\in {\rm \mathbb{N}} \right)\left(\neg \left(\neg \left(\exists k\in {\rm \mathbb{N}} \right)\left(k>1{\rm \; \; }\wedge {\rm \; \; }a(k;n)\in nD_{T} \right){\rm \; \; }\vee {\rm \; \; }\left(\forall r\in {\rm \mathbb{N}} \right)\left(a(r;n)\in nD_{T} \right)\right)\right){\rm \; }\Leftrightarrow \] 
\begin{spacing}{-1.5}\end{spacing}\[\Leftrightarrow {\rm \;}\left(\exists n\in {\rm  \mathbb{N}} \right)\left(\left(\exists k\in {\rm \mathbb{N}} \right)\left(k>1{\rm \; \; }\wedge {\rm \; \; }a(k;n)\in nD_{T} \right){\rm \; \; }\wedge {\rm \; \; }\neg \left(\forall r\in {\rm \mathbb{N}} \right)\left(a(r;n)\in nD_{T} \right)\right){\rm \;}\Leftrightarrow \] 
\begin{spacing}{-1.5}\end{spacing}\[\Leftrightarrow {\rm \;}\left(\exists n\in {\rm  \mathbb{N}} \right)\left(\left(\exists k\in {\rm \mathbb{N}} \right)\left(k>1{\rm \; \; }\wedge {\rm \; \; }a(k;n)\in nD_{T} \right){\rm \; \; }\wedge {\rm \; \; }\left(\exists r\in {\rm \mathbb{N}} \right)\left(a(r;n)\notin nD_{T} \right)\right){\rm \;}\Leftrightarrow \] 
\begin{spacing}{-1.5}\end{spacing}\begin{equation} \label{(16)} 
\Leftrightarrow {\rm \;}\left(\exists n\in {\rm \mathbb{N}} \right)\left(\left(\exists k\in {\rm \mathbb{N}} \right)\left(k>1{\rm \; \; }\wedge {\rm \; \; }a(k;n)\in nD_{T} \right){\rm \; \; }\wedge {\rm \; \; }\left(\exists r\in {\rm \mathbb{N}} \right)\left(a(r;n)\in D_{T} \right)\right).     
\end{equation} 
It follows from \eqref{(16)} that the found element $a(r;n)$ is defining. Let us consider 2 cases.  

\textbf{Case 1.} The element $a(r;n)$ is leading and defining, that is $a(r;n)\in {\rm M} _{T} \backslash \{ 25\} $. Then, by Theorem 4.3, there is an equation $a(r;1)=a(1;n)$.

Using the chain of equivalences from the proof of Theorem 4.3, we get:
\begin{spacing}{0.3}\end{spacing}\begin{equation} \label{(17)} 
p(r)=f(n).                                                                        
\end{equation} 
\begin{equation} \label{(18)} 
k>1{\rm \; \; }\Leftrightarrow {\rm \; \;  }p(k)\ne 5.  
\end{equation} 
\[a(r;n)\in {\rm M} _{T} \backslash \{ 25\} {\rm \; \; }\Rightarrow {\rm \; \; }a(r;n)=p^{2} (r)\ne 25{\rm \; \; }\Rightarrow {\rm \; \; }p(r)\ne 5{\rm \; \; }\mathop{\Leftrightarrow }\limits^{{\rm \eqref{(17)}}} {\rm \; \; }f(n)\ne 5{\rm \; \; }\mathop{\Rightarrow }\limits^{\eqref{(3)},{\rm \; }\eqref{(18)}{\rm \; \; }} \] 
\[\Rightarrow{\rm \; \; } p(k)\cdot f(n)=a(k;n){\rm \; \; }\wedge {\rm \; \; 5}\not {\rm |}a(k;n){\rm \; \; }\mathop{\Leftrightarrow }\limits^{{\rm \eqref{(17)}}} {\rm \; \; }a(k;n)=p(k)\cdot p(r){\rm \; \; }\wedge {\rm \; \; 5}\not {\rm |}a(k;n).\] 

 By Definition 1.1, the latter means that $a(k;n)\in D_{T}$. As a result, a contradiction.

\textbf{Case 2. } The element $a(r;n)$ is not leading, but defining, that is $a(r;n)\in D_{T} \backslash {\rm M} _{T} $.
\begin{spacing}{1.7}\end{spacing}
It follows from Definition 1.2 for element $a(k;n)\mathop{=}\limits^{\eqref{(3)}} p(k)\cdot f(n)\in nD_{T}$, where $k>1$, that
\[\left(\exists p_{1} ,..,p_{s} \in {\rm \mathbb{P}} ,{\rm \; }s\in {\rm \mathbb{N}} \backslash \{ 1\} \right) \left(f(n)=\prod _{j=1}^{s}p_{j}  \right){\rm \; }\vee {\rm \; \; }5|f(n).\] 

 Let us look at these two possible situations.
\begin{spacing}{0.8}\end{spacing}
\[1) \left(\exists p_{1} ,..,p_{s} \in {\rm \mathbb{P}} ,{\rm \; }s\in {\rm \mathbb{N}} \backslash \{ 1\} \right) \left(f(n)=\prod _{j=1}^{s}p_{j}  \right). \]

Note that, by Theorem 3.10, this form for number $f(n)$ is unique. In this situation,
\begin{spacing}{0.8}\end{spacing}\[a(r;n)=p(r)\cdot f(n)=p(r)\cdot \prod _{j=1}^{s}p_{j}  {\rm \; \; }\Rightarrow {\rm \; \; }a(r;n)\in nD_{T} {\rm \; \; }\Leftrightarrow {\rm \; \; }a(r;n)\notin D_{T} .\] 

As a result, a contradiction with the condition $a(r;n)\in D_{T} \backslash {\rm M} _{T} $.
\begin{spacing}{-0.5}\end{spacing}\[2) {\rm \;}5|f(n){\rm \;}\Leftrightarrow {\rm \;}\left(\exists b\in {\rm \mathbb{N}} \right)\left(f(n)=5 \cdot b\right){\rm \;}\Leftrightarrow {\rm \;}\left(\exists b\in {\rm \mathbb{N}} \right)\left(p(r)\cdot f(n)=p(r)\cdot 5 \cdot b\right){\rm \; \; }\mathop{\Leftrightarrow }\limits^{{\rm \eqref{(3)}}} \] 
\begin{spacing}{-1.5}\end{spacing}
\[\Leftrightarrow {\rm \; }\left(\exists b\in {\rm \mathbb{N}} \right)\left(a(r;n)=5\cdot b \cdot p(r)\right){\rm \;}\Rightarrow {\rm \;}5{\rm |}a(r;n){\rm \;}\Rightarrow {\rm \;}a(r;n)\in nD_{T} {\rm \;}\Leftrightarrow {\rm \;}a(r;n)\notin D_{T} .\] 
\begin{spacing}{0.8}\end{spacing}
As a result, a contradiction with the condition $a(r;n)\in D_{T} \backslash {\rm M} _{T} $. 

Thus,
\begin{spacing}{-0.8}\end{spacing}\[\left(\forall n\in {\rm  \mathbb{N}} \right)\left(\left(\exists k\in {\rm  \mathbb{N}} \right)\left(k>1{\rm \; \; }\wedge {\rm \; \; }a(k;n)\in nD_{T} \right){\rm \;}\Rightarrow {\rm \;}\left(\forall r\in {\rm  \mathbb{N}} \right){\rm \; }\left(a(r;n)\in nD_{T} \right)\right).\] 
\begin{spacing}{0.8}\end{spacing}
The necessity is proved.

\textbf{Sufficiency.  }Obviously, \begin{spacing}{-0.5}\end{spacing}\[\left(\forall n\in {\rm \mathbb{N}} \right)\left(\left(\forall r\in {\rm \mathbb{N}} \right)\left(a(r;n)\in nD_{T} \right){\rm \;}\Rightarrow {\rm \;} \left(\exists k\in {\rm \mathbb{N}} \right)\left(k>1{\rm \; \; }\wedge {\rm \; \; }a(k;n)\in nD_{T} \right)\right).\]\begin{spacing}{0.8}\end{spacing}

Theorem 4.4 is proved.

COROLLARY 4.5 (about subcolumns comprising defining elements of $T$-matrix).
\begin{spacing}{0.5}\end{spacing}\[\left(\forall n\in {\rm \mathbb{N}} \right)\left(\left(\exists r\in {\rm \mathbb{N}} \right)\left(a(r;n)\in D_{T} \right){\rm \; }\Leftrightarrow {\rm \;}\left(\forall k\in {\rm \mathbb{N}} \right)\left(k>1{\rm \;}\Rightarrow {\rm \;}a(k;n)\in D_{T} \right)\right).\] \begin{spacing}{0.8}\end{spacing}
\begin{spacing}{1.5}\end{spacing}

PROOF. Use Theorem 4.4. Then, 
\begin{spacing}{-0.5}\end{spacing}\[\left(\forall n\in {\rm \mathbb{N}} \right)\left(\left(\exists k\in {\rm \mathbb{N}} \right)\left(k>1{\rm \; \; }\wedge {\rm \; \; }a(k;n)\in nD_{T} \right){\rm \; }\Leftrightarrow {\rm \; }\left(\forall r\in {\rm \mathbb{N}} \right)\left(a(r;n)\in nD_{T} \right)\right){\rm \; }\Leftrightarrow {\rm \; }\] 
\[\Leftrightarrow {\rm \;}\left(\forall n\in {\rm \mathbb{N}} \right)\left(\neg\left(\forall r\in {\rm \mathbb{N}} \right)\left(a(r;n)\in nD_{T} \right){\rm \; }\Leftrightarrow {\rm \; }\neg \left(\exists k\in {\rm \mathbb{N}} \right)\left(k>1{\rm \; \; }\wedge {\rm \; \; }a(k;n)\in nD_{T} \right)\right){\rm \;}\Leftrightarrow \] 
\[\Leftrightarrow {\rm \;}\left(\forall n\in {\rm \mathbb{N}} \right)\left(\left(\exists r\in {\rm \mathbb{N}} \right)\left(a(r;n)\notin nD_{T} \right){\rm \;}\Leftrightarrow {\rm \;}\left(\forall k\in {\rm \mathbb{N}} \right)\left(\neg \left(k>1{\rm \; \; }\wedge {\rm \; \; }a(k;n)\in nD_{T} \right)\right)\right){\rm \;}\Leftrightarrow \] 
\[\Leftrightarrow {\rm \;}\left(\forall n\in {\rm \mathbb{N}} \right)\left(\left(\exists r\in {\rm \mathbb{N}} \right)\left(a(r;n)\notin nD_{T} \right){\rm \;}\Leftrightarrow {\rm \;}\left(\forall k\in {\rm \mathbb{N}} \right)\left(\neg \left(k>1\right){\rm \; \; }\vee {\rm \; \; }a(k;n)\notin nD_{T} \right)\right){\rm \;}\Leftrightarrow \] 
\begin{spacing}{-1.2}\end{spacing}
\[\Leftrightarrow {\rm \;}\left(\forall n\in {\rm \mathbb{N}} \right)\left(\left(\exists r\in {\rm \mathbb{N}} \right)\left(a(r;n)\notin nD_{T} \right){\rm \;}\Leftrightarrow {\rm \; }\left(\forall k\in {\rm \mathbb{N}} \right)\left(k>1{\rm \;}\Rightarrow {\rm \; }a(k;n)\notin nD_{T} \right)\right){\rm \; }\Leftrightarrow \] 
\begin{spacing}{-1.5}\end{spacing}
\[\Leftrightarrow {\rm \;}\left(\forall n\in {\rm \mathbb{N}} \right)\left(\left(\exists r\in {\rm \mathbb{N}} \right)\left(a(r;n)\in D_{T} \right){\rm \;}\Leftrightarrow {\rm \;}\left(\forall k\in {\rm \mathbb{N}} \right)\left(k>1{\rm \;}\Rightarrow {\rm \; }a(k;n)\in D_{T} \right)\right).\] \begin{spacing}{0.8}\end{spacing}

Corollary 4.5 is proved.

LEMMA 4.6.
\begin{spacing}{0.3}\end{spacing}
\[\left(\forall n\in {\rm \mathbb{N}} \right)\left(\left(\exists r\in {\rm \mathbb{N}} \right)\left(a(r;n)\in D_{T} \right){\rm \;}\Rightarrow {\rm \;}\left(\exists !k\in {\rm \mathbb{N}} \right)\left(a(k;n)\in {\rm M} _{T} \right)\right).\]

PROOF. Choose any $n$-column of $T$-matrix, where $n\in{\rm \mathbb{N}}$, in which condition of Lemma 4.6 is satisfied. Then from Corollary 4.5 we get:
\begin{spacing}{0.5}\end{spacing}
\[\left(\forall k\in {\rm \mathbb{N}} \backslash \{ 1\} \right)\left(a(k;n)\in D_{T} \right).\]

Given Definition 1.1, this means: 
\begin{spacing}{0.3}\end{spacing}\[\left(\forall k\in {\rm \mathbb{N}} \backslash \{ 1\} \right)\left(5\not |a(k;n){\rm \; \; }\wedge {\rm \; \; }\left(\exists p_{1} ,p_{2} \in {\rm  \mathbb{P}} \right)\left(a(k;n)=p_{1} \cdot p_{2} \right)\right).\] \begin{spacing}{0.8}\end{spacing}

 It follows that 
\begin{spacing}{0.1}\end{spacing}
 \[a(k;n)=p(k)\cdot f(n){\rm \; \;}\wedge {\rm \; \;}p(k),f(n)\in {\rm  \mathbb{P}} \backslash \{ 2;3;5\} .\]
 
Now if to go over all $k$-elements of $n$-column of $T$-matrix, where $k\in {\rm \mathbb{N}} \backslash \{ 1\}$, then given Lemma 1.1 and Property 3.1,
\begin{spacing}{0.3}\end{spacing}\[\left(\exists !k\in {\rm \mathbb{N}} \right)\left(p(k)=f(n)\right){\rm \;}\Leftrightarrow {\rm \;}\left(\exists !k\in {\rm \mathbb{N}} \right)\left(a(k;n)=p^{2} (k)\right).\] \begin{spacing}{1}\end{spacing}

By Definition 1.3, the latter means that
\begin{spacing}{0.3}\end{spacing}\[\left(\exists !k\in {\rm \mathbb{N}} \right)\left(a(k;n)\in {\rm M} _{T} \right).\] \begin{spacing}{0.8}\end{spacing}

Lemma 4.6 is proved.

LEMMA 4.7.
\begin{spacing}{0.2}\end{spacing}
\[\left(\forall n\in {\rm \mathbb{N}} \backslash \{ 1\} \right)\left(\left(\exists r\in {\rm \mathbb{N}} \right)\left(a(r;n)\in D_{T} \right){\rm \;}\Leftrightarrow {\rm \;}f(n)\in {\rm \mathbb{P}} \right).\]

PROOF. Choose any $n$-column of $T$-matrix, where $n\in {\rm \mathbb{N}} \backslash \{ 1\} $. 

\textbf{Necessity. } Let $\left(\exists r\in {\rm \mathbb{N}} \right)\left(a(r;n)\in D_{T} \right)$. Then it follows from Definition 1.1 for element $a(r;n)$ that
\begin{spacing}{0}\end{spacing}\[\left(\exists r\in {\rm \mathbb{N}} \right)\left(a(r;n)=p(r)\cdot f(n){\rm \; }\wedge {\rm \;}p(r),{\rm \; }f(n)\in {\rm \mathbb{P}} \backslash \{ 2;3;5\} \right){\rm \; }\Rightarrow\]
\[\Rightarrow {\rm \; }f(n)\in {\rm \mathbb{P}} \backslash \{ 2;3;5\} {\rm \;}\mathop{\Rightarrow }\limits^{{\rm \mathbb{P}} \backslash \{ 2;3;5\} \subset {\rm \mathbb{P}}} f(n)\in {\rm \mathbb{P}} .\] 

The necessity is proved. 

\textbf{Sufficiency. } Let $f(n)\in {\rm  \mathbb{P}}$, where $n\in {\rm \mathbb{N}}\backslash \{ 1\}$. It follows from Lemma 1.1 that
\begin{spacing}{0.5}\end{spacing}
\[f(n)>5, n\in {\rm \mathbb{N}}\backslash \{ 1\}.\]

 In turn, it follows from rule \eqref{(2)} that ${\rm \; }p(r)\in {\rm \mathbb{P}} \backslash \{2;3;5\}, r\in {\rm \mathbb{N}} \backslash \{ 1\}$.

 Then given Definition 1.1 for some number $r\in {\rm \mathbb{N}} \backslash \{ 1\} $:
\begin{spacing}{0.5}\end{spacing}
\[p(r)\cdot f(n)\mathop{=}\limits^{\eqref{(3)}} a(r;n)\in D_{T} .\] 
\[\left(\exists r\in {\rm \mathbb{N}} \backslash \{ 1\} \right)\left(a(r;n)\in D_{T} \right){\rm \; \; }\mathop{\Rightarrow }\limits^{{\rm \mathbb{N}} \backslash \{ 1\} \subset {\rm \mathbb{N}} } {\rm \; \; }\left(\exists r\in {\rm \mathbb{N}} \right)\left(a(r;n)\in D_{T} \right).\] 

Lemma 4.7 is proved. 

CONCLUSION 4.1. It follows from Corollary 4.5 and Lemma 4.7 that each subcolumn comprising defining elements of $T$-matrix, that is each column comprising defining elements of matrix $T^*$, one-to-one corresponds with some prime number.

CONCLUSION 4.2. Each column comprising not defining elements of $T$-matrix one-to-one corresponds with some composite number of the form $6h\pm 1,{\rm \; }h\in {\rm\mathbb{N}} $. 

Thus, several theorems about prime numbers can be formulated in terms of $T$-matrix elements.

THEOREM 4.8. Sets $\widetilde{T},{\rm \; }{\rm M} _{T} ,{\rm \; }nD_{T} ,{\rm \; }D_{T} $ are countable.

PROOF. 1) By Property 2.1, $T$-matrix is matrix of infinite order. So, $T$-matrix has countably infinite number of rows and columns. On the basis of Theorem 3.6, we get that set $\widetilde{T}$ is countable.

2) By Lemma 4.2, each row of $T$-matrix contains a unique leading element. It follows from Property 2.1 that $T$-matrix has countably infinite number of rows. Then set ${\rm M} _{T} $ of all leading elements of $T$-matrix is countable.

3) It follows from Property 2.4 that all elements of the first row of $T$-matrix are not defining. By Property 2.1, $T$-matrix is matrix of infinite order. Therefore, any row of $T$-matrix contains countably infinite number of elements. In particular, set $nD_{T_{1} } $ is countable. Given Theorem 4.8, 1) and Theorem 3.8, clear that set $nD_{T} \backslash nD_{T_{1} } \subset \widetilde{T}$ of not defining elements of $T$-matrix, which are not included in set $nD_{T_{1} } $, is not more than countable. Then by Theorem 3.9, set 
\begin{spacing}{0.5}\end{spacing}
\[nD_{T} =nD_{T_{1} } \cup {\rm \; } (nD_{T} \backslash nD_{T_{1} } )\]

\noindent is countable.

4) The first row of $T$-matrix doesn't contain the defining elements. As such, consider the matrix $T^*$. Using Theorem 4.4 and Corollary 4.5, we separate the columns comprising defining elements from the columns comprising not defining elements in matrix  $T^*$. It follows from Theorem 3.4 and Conclusion 4.1 that there is a countable number of columns comprising defining elements in matrix $T^*$. It follows from Property 2.1 that $T$-matrix has countably infinite number of rows and columns. Then, by Theorem 3.5, columns comprising defining elements in matrix $T^*$ contain countably infinite number of elements. Using Theorem 3.6, in result we get that set $D_{T} $ is countable. 

Theorem 4.8 is proved.

THEOREM 4.9. Sets $D_{T,6h+1} $, $D_{T,6h-1} $ are countable.

PROOF. It follows from Theorem 3.2, Theorem 3.3 and Conclusion 4.1 that there is a countable number of columns comprising defining elements in matrix $T^*$ and one-to-one corresponding to some prime numbers $f(n)=6h_{{\rm \; }1} \pm 1$. In turn, by Lemma 4.1, each such column contains the elements of the form $6h\pm 1$. Then given the proof of Lemma 4.1, there are equivalences with respect to $f(n)=6h_{{\rm \; }1} \pm 1$:
\begin{spacing}{0.3}\end{spacing}\[1){\rm \;} p(k)=6t+1{\rm \; \; }\Leftrightarrow {\rm \; \; }a(k;n)=6h\pm 1.\] 
\begin{spacing}{-0.3}\end{spacing}\[2){\rm \;} p(k)=6t-1{\rm \; \; }\Leftrightarrow {\rm \; \; }a(k;n)=6h\mp 1.\] \begin{spacing}{1.2}\end{spacing}

It follows from Theorem 3.2 and Theorem 3.3 that the set of elements of the form $6h+1${\rm \; }$(h\in {\rm \mathbb{N}}) $ and the set of elements of the form $6h-1 {\rm \; }(h\in {\rm \mathbb{N}})$ from each column under review are countable. Then, by Theorem 3.6, sets $D_{T,6h+1}$, $D_{T,6h-1}$ are countable. 

Theorem 4.9 is proved.

COMMENT. $D_{T,6h+1} \cup D_{T,6h-1} =D_{T} $. Then, by Theorem 3.7, set $D_{T} $ is countable.

THEOREM 4.10. Sets $nD_{T,6h+1} $, $nD_{T,6h-1} $ are countable.

PROOF. Note that the first row of $T$-matrix contains only not defining elements. From the proof of Lemma 1.1, it is clear that
\begin{spacing}{0.5}\end{spacing}
\[1) {\rm \; }n=2h_{1} -1, h_{1} \in {\rm \mathbb{N}} {\rm \;}\Leftrightarrow {\rm \;}f(n)=6h_{1} -1{\rm \;}\Leftrightarrow \]
\[ \Leftrightarrow{\rm \;}p(1)\cdot f(n)=p(1)\cdot (6h_{1} -1){\rm \; \;}\mathop{\Leftrightarrow }\limits^{\eqref{(3)}} {\rm \; \; }a(1;n)=5\cdot (6h_{1} -1).\] 
\[5\cdot (6h_{1} -1)=5\cdot (6h_{1} )-5=6\cdot (5h_{1} )-6+1=\]
\[=6\cdot (5h_{1} -1)+1=6h+1, h=5h_{1} -1.\]
\[n=2h_{1} -1, h_{1} \in {\rm \mathbb{N}} {\rm \; }\Leftrightarrow {\rm \;}a(1;n)=6h+1.\] 
\[2){\rm \; } n=2h_{2} , h_{2} \in {\rm \mathbb{N}} {\rm \; }\Leftrightarrow {\rm \; }f(n)=6h_{2} +1{\rm \;}\Leftrightarrow \]
\[\Leftrightarrow{\rm \;}p(1)\cdot f(n)=p(1)\cdot (6h_{2} +1){\rm \; \;}\mathop{\Leftrightarrow }\limits^{\eqref{(3)}} {\rm \; \;}a(1;n)=5\cdot (6h_{2} +1).\] 
\[5\cdot (6h_{2} +1)=5\cdot (6h_{2} )+5=6\cdot (5h_{2} )+6-1=\]
\[=6\cdot (5h_{2} +1)-1=6h-1, h=5h_{2} +1.\]
\[n=2h_{2} , h_{2} \in {\rm \mathbb{N}} {\rm \;}\Leftrightarrow {\rm \;}a(1;n)=6h-1.\] 
\begin{spacing}{1.2}\end{spacing}
It follows that set $nD_{T_{1} ,6h+1}$ of all $T$-matrix not defining elements of the first row and the form $6h+1$ is equinumerous with the set of all odd numbers, set $nD_{T_{1} ,6h-1} $ of all $T$-matrix not defining elements of the first row and the form $6h-1$ is equinumerous with the set of all even numbers. The set of all odd numbers and the set of all even numbers are countable. Then, sets $nD_{T_{1} ,6h+1} $, $nD_{T_{1} ,6h-1} $ are also countable. By Theorem 4.8, 3), set $nD_{T} $ is countable. By Theorem 3.8, sets  $nD_{T,6h+1} ,{\rm \; }nD_{T,6h-1} \subset nD_{T} $ are not more than countable. Then, on the basis of Theorem 3.9, sets 
\begin{spacing}{0.5}\end{spacing}\[nD_{T,6h+1} =nD_{T,6h+1} \cup nD_{T_{1} ,6h+1} , nD_{T,6h-1} =nD_{T,6h-1} \cup nD_{T_{1} ,6h-1} \] \begin{spacing}{1.2}\end{spacing}
\noindent are countable.

Theorem 4.10 is proved.

THEOREM 4.11 (about the «transition down» of $T$-matrix defining element). 

Let $\# _{k} \left(a\right)$ is a number of element $a$ in $k$-row of $T$-matrix. Then
\begin{spacing}{0.5}\end{spacing}\[\left(\forall k;n\in {\rm \mathbb{N}} \right)\left(p^{2} (k)<a(k;n){\rm \; \; }\wedge {\rm \; \; }a(k;n)\in D_{T} {\rm \;}\Rightarrow \right. \] 
\[\left. \Rightarrow {\rm \;}\left(\exists {\rm \; }j\in {\rm \mathbb{N}}\right)\left(k<j{\rm \; \;}\wedge {\rm \; \; }a(k;n)<p^{2} (j){\rm \; \; }\wedge {\rm \; \; }a(j;\# _{k} (p^{2} (k)))=a(k;n){\rm \; \; }\wedge {\rm \; \;}a(j;n)=p^{2} (j)\right)\right).\] \begin{spacing}{1.2}\end{spacing}

PROOF. Choose any defining element $a(k;n)$ that is larger than leading element $p^{2} (k)$ in $k$-row of $T$-matrix. Then it follows from Definition 1.1 that
\begin{spacing}{0.2}\end{spacing}\[a(k;n)=p(k)\cdot f(n){\rm \; \;}\wedge {\rm \; \;}p(k),{\rm \; }f(n)\in {\rm \mathbb{P}} \backslash \{ 2;3;5\} .\] 

\begin{spacing}{1}\end{spacing}
Therefore, prime number $f(n)>5$ is an element of sequence $\left(p(k)\right)_{k=1}^{\infty } $. As such,
\begin{spacing}{0.4}\end{spacing}
\begin{equation} \label{(19)} 
\left(\exists {\rm \; }j\in {\rm \mathbb{N}} \right)\left(p(j)=f(n)\right).  
\end{equation} 
\begin{spacing}{1.7}\end{spacing}
$1) {\rm \;} p^{2} (k)<a(k;n){\rm \; \; }\mathop{\Leftrightarrow }\limits^{{\rm \eqref{(3)},\; \eqref{(19)}}} {\rm \; \; }p^{2} (k)<p(k)\cdot p(j){\rm \; \; }\Leftrightarrow {\rm \; \; }p(k)<p(j){\rm \; \; }\Leftrightarrow {\rm \; \; }k<j.$
\begin{spacing}{2}\end{spacing}
$2)  {\rm \;} p(k)<p(j){\rm \; \; }\Leftrightarrow {\rm \; \; }p(k)\cdot p(j)<p^{2} (j){\rm \; \;}\mathop{\Leftrightarrow }\limits^{{\rm \eqref{(19)},\;  \eqref{(3)}}} {\rm \; \; \; }a(k;n)<p^{2} (j).$
\begin{spacing}{2}\end{spacing}
$3) {\rm \;}  a(j;\# _{k} (p^{2} (k)))\mathop{=}\limits^{ \eqref{(3)}} p(j)\cdot f(\# _{k} (p^{2} (k)))=p(j)\cdot p(k)=p(k)\cdot p(j)\mathop{=}\limits^{\eqref{(19)}} p(k)\cdot f(n)\mathop{=}\limits^{ \eqref{(3)}} a(k;n){\rm \;}\Rightarrow$
\begin{spacing}{0.5}\end{spacing}\[\Rightarrow {\rm \; \; }a(j;\# _{k} (p^{2} (k)))=a(k;n).\]\begin{spacing}{0.8}\end{spacing}
\begin{spacing}{2}\end{spacing}
$4)  {\rm \;} a(j;n)\mathop{=}\limits^{ \eqref{(3)}} p(j)\cdot f(n){\rm \; \; }\mathop{=}\limits^{\eqref{(19)}} {\rm \; }p^{2} (j){\rm \; \; }\Rightarrow {\rm \; \; }a(j;n)=p^{2} (j).$

Theorem 4.11 is proved.

THEOREM 4.12. If between any two columns of matrix $T^*$  -- $A$ and $B$, consisting of not defining elements, there are $m>1$ consecutive columns $C_{1} ,...,C_{m} $, consisting of defining elements, then between columns $A$ and $B$ there exists a singular submatrix $T'$ of finite order $m$ of matrix $T^*$, and the leading elements lie on the main diagonal of matrix $T'$.

PROOF. Choose any such two columns $A$ and $B$. The proof is by induction on number $m$ of consecutive columns $C_{1} ,...,C_{m} $.

\textbf{Base case. }$m=2$. By Lemma 4.6, column $C_{1} $ contains a unique leading element $p^{2} (k)$ for some number $k\in {\rm \mathbb{N}} $. Since $C_{1} ,{\rm \; }C_{2} $ are the columns consisting of defining elements, from Theorem 4.11 and Conclusion 4.1 we get a matrix $T'$ in the form:

\[T'=\mathop{\left(\begin{array}{cc} {p^{2} (k)} & {p(k)\cdot p(k+1)} \\ {p(k+1)\cdot p(k)} & {p^{2} (k+1)} \end{array}\right)}\limits^{C_{1} {\rm \; \; \; \; \; \; \; \; \; \; \; \; \; \; \; \; \; \; \; \; \; \; \; \; \; \; \; \; \; \; \;}C_{2} } .\] 

 Thus, leading elements $p^{2} (k)$, $p^{2} (k+1)$ lie on the main diagonal of matrix $T'$. Now we find the determinant of matrix $T'$. 
\begin{spacing}{0.7}\end{spacing}\[\det T'=\left|\begin{array}{cc} {p^{2} (k)} & {p(k)\cdot p(k+1)} \\ {p(k+1)\cdot p(k)} & {p^{2} (k+1)} \end{array}\right|=\]
\[=p^{2} (k)\cdot p^{2} (k+1)-p(k+1)\cdot p(k)\cdot p(k)\cdot p(k+1)=0.\] \begin{spacing}{1.0}\end{spacing}

As a result, the square matrix $T'$ of finite order 2 is singular.

\textbf{Inductive hypothesis. } Suppose that the statement of Theorem 4.12 is true for any number $m\in {\rm \mathbb{N}} \backslash \{ 1\} $.

\textbf{Inductive step. }Show that the statement of Theorem 4.12 is true for $m+1$ consecutive columns $C_{1} ,...,C_{m+1} $. 

Separate the submatrix $T'$ from matrix $T^*$:
\[\mathop{\left(\begin{array}{ccccc} {p^{2} (k)} & {a(k;n+1)} & {...} & {a(k;n+m-1)} & {a(k;n+m)} \\ {a(k+1;n)} & {p^{2} (k+1)} & {...} & {a(k+1;n+m-1)} & {a(k+1;n+m)} \\ {a(k+2;n)} & {a(k+2;n+1)} & {...} & {a(k+2;n+m-1)} & {a(k+2;n+m)} \\ {...} & {...} & {...} & {...} & {...} \\ {a(k+m-1;n)} & {a(k+m-1;n+1)} & {...} & {p^{2} (k+m-1)} & {a(k+m-1;n+m)} \\ {a(k+m;n)} & {a(k+m;n+1)} & {...} & {a(k+m;n+m-1)} & {a(k+m;n+m)} \end{array}\right)}\limits^{\mathop{C_{1} }\limits_{\downarrow } {\rm \; \; \; \; \; \; \; \; \; \; \; \; \; \; \; \; \; \; \; \; \; \; \; \; \; \; \; \; \; \; \; \; \; \; }\mathop{C_{2} }\limits_{\downarrow }{\rm \; \; \; \; \; \; \; \; \; \; \; \; \; \; \; \; \; \; \; \; \; \; \; \; \; \; \; \; \; \; \; \; \; \; \; \; \; \; \; \; \; \; \; \; \; \; \; \; \; \; \; \; \; \; \; \; }\mathop{C_{m} }\limits_{\downarrow } {\rm \; \; \; \; \; \; \; \; \; \; \; \; \; \; \; \; \; \; \; \; \; \; \; \; \; \; \; \; \; \; \; \; \; \; \; \; }\mathop{C_{m+1} }\limits_{\downarrow } } .\]

By the induction hypothesis, the submatrix $T''$ of matrix $T^*$ (of matrix $T'$):
\[\mathop{\left(\begin{array}{ccccc} {p^{2} (k)} & {a(k;n+1)} & {a(k;n+2)} & {...} & {a(k;n+m-1)} \\ {a(k+1;n)} & {p^{2} (k+1)} & {a(k+1;n+2)} & {...} & {a(k+1;n+m-1)} \\ {a(k+2;n)} & {a(k+2;n+1)} & {p^{2} (k+2)} & {...} & {a(k+2;n+m-1)} \\ {...} & {...} & {...} & {...} & {...} \\ {a(k+m-1;n)} & {a(k+m-1;n+1)} & {a(k+m-1;n+2)} & {...} & {p^{2} (k+m-1)} \end{array}\right)}\limits^{{\rm \; \; \; \; \; \; \; \; \; \; \; \; }\mathop{C_{1} }\limits_{\downarrow } {\rm \; \; \; \; \; \; \; \; \; \; \; \; \; \; \; \; \; \; \; \; \; \; \; \; \; \; \; \; \; \; \; \; \; \; \;  }\mathop{C_{2} }\limits_{\downarrow } {\rm \; \; \; \; \; \; \; \; \; \; \; \; \; \; \; \; \; \; \; \; \; \; \; \; \; \; \; \; \; \; \; \; \; \; \; \; \; \; \; \; \; \; \; \; }\mathop{C_{3} }\limits_{\downarrow } {\rm \; \; \; \; \; \; \; \; \; \; \; \; \; \; \; \; \; \; \; \; \; \; \; \; \; \; \; \; \; \; \; \; \; \; \; \; \; \; \; \; \; \; \; \; \; \; \; \; \;}\mathop{C_{m} }\limits_{\downarrow } {\rm \; \; \; \; \; \; \; \; \; \; \; \; \; \; }} ,\] 
is singular matrix of finite order $m$, and the leading elements $p^{2} (k+i),{\rm \; }i=\overline{0;m-1}$ lie on the main diagonal of matrix $T''$.

 1) Use Theorem 4.11 and the presentation of matrix $T'$:
\begin{spacing}{0.5}\end{spacing}\[p^{2} (k+m-1)<a(k+m-1;n+m){\rm \; \; }\wedge {\rm \; \; }a(k+m-1;n+m)\in D_{T} {\rm \; \; }\Rightarrow \] 
\[\Rightarrow {\rm \; }k+m-1<k+m{\rm \; \; \; }\wedge {\rm \; \; \; }a(k+m-1;n+m)<p^{2} (k+m){\rm \; \; \; }\wedge \] 
\[\wedge {\rm \; \; \; }a(k+m;{\rm \; }\# _{k+m-1} (p^{2} (k+m-1)))=a(k+m-1;n+m){\rm \; \; }\wedge \]
\begin{equation} \label{(20)} 
\wedge {\rm \; \; \; }a(k+m;n+m)=p^{2} (k+m){\rm \; \; }\Rightarrow {\rm \; \; }a(k+m;n+m)=p^{2} (k+m).
\end{equation} \begin{spacing}{1.1}\end{spacing}

2) Now we find the determinant of matrix $T'$. Note that given \eqref{(20)},
\begin{spacing}{0.5}\end{spacing}\[p^{2} (k+i)=a(k+i;n+i),{\rm \; }i=\overline{0;m}.\] 
\[\det T'=\left|\begin{array}{ccccc} {a(k;n)} & {a(k;n+1)} & {a(k;n+2)} & {...} & {a(k;n+m)} \\ {a(k+1;n)} & {a(k+1;n+1)} & {a(k+1;n+2)} & {...} & {a(k+1;n+m)} \\ {a(k+2;n)} & {a(k+2;n+1)} & {a(k+2;n+2)} & {...} & {a(k+2;n+m)} \\ {...} & {...} & {...} & {...} & {...} \\ {a(k+m;n)} & {a(k+m;n+1)} & {a(k+m;n+2)} & {...} & {a(k+m;n+m)} \end{array}\right|\mathop{=}\limits^{\eqref{(3)}} \] 
\begin{spacing}{0.5}\end{spacing}
\[=\left|\begin{array}{ccccc} {p(k)\cdot f(n)} & {p(k)\cdot f(n+1)} & {a(k;n+2)} & {...} & {a(k;n+m)} \\ {p(k+1)\cdot f(n)} & {p(k+1)\cdot f(n+1)} & {a(k+1;n+2)} & {...} & {a(k+1;n+m)} \\ {p(k+2)\cdot f(n)} & {p(k+2)\cdot f(n+1)} & {a(k+2;n+2)} & {...} & {a(k+2;n+m)} \\ {...} & {...} & {...} & {...} & {...} \\ {p(k+m)\cdot f(n)} & {p(k+m)\cdot f(n+1)} & {a(k+m;n+2)} & {...} & {a(k+m;n+m)} \end{array}\right|=\] 
\begin{spacing}{0.3}\end{spacing}
\[=f(n)\cdot f(n+1)\cdot \left|\begin{array}{ccccc} {p(k)} & {p(k)} & {a(k;n+2)} & {...} & {a(k;n+m)} \\ {p(k+1)} & {p(k+1)} & {a(k+1;n+2)} & {...} & {a(k+1;n+m)} \\ {p(k+2)} & {p(k+2)} & {a(k+2;n+2)} & {...} & {a(k+2;n+m)} \\ {...} & {...} & {...} & {...} & {...} \\ {p(k+m)} & {p(k+m)} & {a(k+m;n+2)} & {...} & {a(k+m;n+m)} \end{array}\right|=\] 
\begin{spacing}{0.5}\end{spacing}
\begin{equation} \label{(21)} 
=f(n)\cdot f(n+1)\cdot 0=0{\rm \; \; }\Rightarrow {\rm \; \; }\det T'=0.  
\end{equation} 

It follows from \eqref{(20)} and \eqref{(21)} that matrix $T'$ is singular matrix of finite order $m+1$, and the leading elements $p^{2} (k+i),{\rm \; }i=\overline{0;m}$ lie on the main diagonal of matrix $T'$. 

Theorem 4.12 is proved.

COMMENT. There exist 2 columns $Q_{1}$ and $Q_{2}$, consisting of not defining elements of matrix $T^*$, such that the columns, consisting of defining elements of matrix $T^*$, don't exist between them. 

It is clear that there exist 2 composite numbers $q_{1}, q_{2}$ of the form $6h\pm 1$ such that primes don't exist between them. Then, on the basis of Theorem 4.2, there is a conjunction:  
\begin{spacing}{0.5}\end{spacing}
\[Q_{1} \leftrightarrow q_{1} {\rm \; \; }\wedge {\rm \; \; }Q_{2} \leftrightarrow q_{2}.\]
So, given Conclusion 4.1, we conclude that the columns, consisting of defining elements of matrix $T^*$, don't exist between columns $Q_{1}$ and $Q_{2}$.

PROPOSITION 4.13. Any real number $x$ can be uniquely expressed as a sum of integer part (entire) and fractional part (mantissa) of number $x$ (see [1]):
\begin{spacing}{0.1}\end{spacing}\[x=\left\lfloor x\right\rfloor +\{ x\} . \]\begin{spacing}{0.8}\end{spacing}

\textbf{}

PROPERTY 4.14 (of number's entire).
\begin{spacing}{0.1}\end{spacing}\begin{equation} \label{(22)}
\left(\forall x\in {\rm \mathbb{R}}\right)\left(\left\lfloor x\right\rfloor \le x\right).
\end{equation}\begin{spacing}{0.8}\end{spacing}

PROPERTY 4.15 (of number's mantissa).
\begin{spacing}{0.1}\end{spacing}
\[\left(\forall x\in {\rm \mathbb{R}}\right)\left(0\le \{ x\} <1{\rm \; }\right).\]
\begin{spacing}{0.8}\end{spacing}

PROPERTY 4.16. 
\begin{spacing}{0.1}\end{spacing} \begin{equation} \label{(23)}
\left(\forall x\in {\rm \mathbb{R}: } {\rm \;} x\ge 0\right)\left(\pi (x)=\pi \left(\left\lfloor x\right\rfloor \right)\right) \text { (see [6])}.
\end{equation}

PROPERTY 4.17.
\begin{spacing}{0.1}\end{spacing} \begin{equation} \label{(24)}
 \left(\forall x\in {\rm \mathbb{R}: }{\rm \;}x\ge 0\right) \left(\pi _{{\rm M} _{T} } (x)=\pi _{{\rm M} _{T} } \left(\left\lfloor x\right\rfloor \right)\right).
\end{equation}

PROOF. Choose any real number $x\ge 0$. Then, using Proposition 4.13 and Property 4.15, we get:
\begin{spacing}{0.1}\end{spacing}\[\pi _{{\rm M} _{T} } (x)=\pi _{{\rm M} _{T} } \left(\left\lfloor x\right\rfloor +\{ x\} \right)=\pi _{{\rm M} _{T} } \left(\left\lfloor x\right\rfloor \right).\] 

Property 4.17 is proved.

Consider now theorems establishing the relationship between $\pi (x)$ and $\pi _{{\rm M} _{T} } (x)$.

THEOREM 4.18.
\begin{spacing}{0.1}\end{spacing} \begin{equation} \label{(25)}
\left(\forall x\in {\rm \mathbb{R}}\right)\left(x\ge 3{\rm \;}\Rightarrow {\rm \;}\pi (x)=\pi _{{\rm M} _{T} } \left(x^{2} \right)+2\right).         
\end{equation}

PROOF. Choose any number $x\in {\rm \mathbb{R}}: x\ge 3$. Then, by Property 4.17,
\begin{spacing}{0.3}\end{spacing}\[\pi _{{\rm M} _{T} } (x^{2} )=\pi _{{\rm M} _{T} } \left(\left\lfloor x^{2} \right\rfloor \right). \]

By Lemma 4.2, each $k$-row of $T$-matrix contains a unique leading element of the form $p^{2} (k)$. Choose any $k$-row of $T$-matrix such that $p^{2} (k)\le \left\lfloor x^{2} \right\rfloor $. Then there is a chain of equivalences:
\begin{spacing}{0.5}\end{spacing}\[p^{2} (k)\le \left\lfloor x^{2} \right\rfloor {\rm \; \; }\Leftrightarrow {\rm \; \; }p^{2} (k)\le x^{2} {\rm \; \; }\Leftrightarrow {\rm \; \; }p(k)\le x{\rm \; \; }\Leftrightarrow {\rm \; \; }p(k)\le \left\lfloor x\right\rfloor .\] 

 Thus, $p^{2} (k)\le \left\lfloor x^{2} \right\rfloor {\rm \;}\Leftrightarrow {\rm \;}p(k)\le \left\lfloor x\right\rfloor $. In turn, the number $\left\lfloor x\right\rfloor$ is greater than or equal to prime numbers 2 and 3. As a result,
\begin{spacing}{0.5}\end{spacing}\[\pi \left(\left\lfloor x\right\rfloor \right)=\pi _{{\rm M} _{T} } \left(\left\lfloor x^{2} \right\rfloor \right)+2{\rm \; \; }\mathop{\Leftrightarrow }\limits^{{\rm \eqref{(23)},\; \eqref{(24)}}} {\rm \; \; }\pi (x)=\pi _{{\rm M} _{T} } \left(x^{2} \right)+2.\]

Theorem 4.18 is proved. 

THEOREM 4.19.
\begin{spacing}{0.5}\end{spacing}\[\left(\forall x\in {\rm \mathbb{R}}\right)\left(x^{2} \ge 3{\rm \; }\Rightarrow {\rm \;}\pi \left((x+1)^{2} \right)-\pi \left(x^{2} \right)=\pi _{{\rm M} _{T} } \left((x+1)^{4} \right)-\pi _{{\rm M} _{T} } \left(x^{4} \right)\right).\]

PROOF. Choose any number $x\in {\rm \mathbb{R}} : x^{2} \ge 3$. Note that $\pi ((x+1)^{2} )-\pi (x^{2} )$ is a number of primes between $x^{2} $ and $(x+1)^{2} $. Therefore, $\pi _{{\rm M} _{T} } ((x+1)^{4} )-\pi _{{\rm M} _{T} } (x^{4} )$ is a number of leading elements of $T$-matrix between $x^{4}$ and $(x+1)^{4}$. Using Theorem 4.18, we get:
\begin{spacing}{0.5}\end{spacing}\[\pi \left((x+1)^{2} \right)-\pi \left(x^{2} \right){\rm \; }\mathop{=}\limits^{(x+1)^{2} >3} {\rm \;}\left(\pi _{{\rm M} _{T} } \left((x+1)^{4} \right)+2\right)-\left(\pi _{{\rm M} _{T} } \left(x^{4} \right)+2\right)=\] 
\[=\pi _{{\rm M} _{T} } \left((x+1)^{4} \right)+2-\pi _{{\rm M} _{T} } \left(x^{4} \right)-2=\pi _{{\rm M} _{T} } \left((x+1)^{4} \right)-\pi _{{\rm M} _{T} } \left(x^{4} \right). \] 

Theorem 4.19 is proved.

PROPOSITION 4.20. There are 2 constants $a$ and $A$ such that $0<a<1,{\rm \; }A>1$, and
\begin{spacing}{0.8}\end{spacing}
\[2 \cdot a \cdot \frac{\sqrt{x} }{\ln x}-2<\pi _{{\rm M} _{T} } (x)<2 \cdot A\cdot \frac{\sqrt{x} }{\ln x}-2 {\rm \;}\text{ for all } x\ge 9.\] 
\begin{spacing}{1.5}\end{spacing}
PROOF. It is known that there are 2 constants $a$ and $A$ such that $0<a<1,{\rm \; }A>1$, and the Chebyshev inequalities (see [4]) are satisfied:
\begin{spacing}{0.6}\end{spacing}\begin{equation} \label{(26)} 
a\cdot \frac{x}{\ln x} <\pi (x)<A\cdot \frac{x}{\ln x}{\rm \; } \text{ for all }x\ge 2.       
\end{equation} \begin{spacing}{1.5}\end{spacing}
Therefore, these inequalities are true for all $x\ge 3$. Then, given Theorem 4.18, 
\[a\cdot \frac{x}{\ln x} <\pi _{{\rm M} _{T} } \left(x^{2} \right)+2<A\cdot \frac{x}{\ln x} ,{\rm \; \; }x\ge 3{\rm \; \; }\Leftrightarrow \]
\[\Leftrightarrow{\rm \; \; }a\cdot \frac{x}{\ln x}-2<\pi _{{\rm M} _{T} } \left(x^{2} \right)<A\cdot \frac{x}{\ln x}-2,{\rm \; }x\ge 3{\rm \; \; }\Leftrightarrow \] 
\[\Leftrightarrow {\rm \; \; }a\cdot \frac{\sqrt{x} }{\ln \sqrt{x} }-2<\pi _{{\rm M} _{T} } (x)<A\cdot \frac{\sqrt{x} }{\ln \sqrt{x} }-2,{\rm \; }x\ge 9{\rm \; \; }\Leftrightarrow  \] 
\[\Leftrightarrow {\rm \; \; }2 \cdot a \cdot \frac{\sqrt{x}}{\ln x}-2<\pi _{{\rm M} _{T} } (x)<2 \cdot A \cdot \frac{\sqrt{x}}{\ln x}-2,{\rm \; }x\ge 9.\]
\begin{spacing}{1.5}\end{spacing}

Proposition 4.20 is proved.

COROLLARY 4.21. There are 2 positive constants $b$ and $B$ such that 
\begin{spacing}{0.5}\end{spacing}
\[b\cdot (k+2)^{2} \cdot \ln ^{2} (k+2)<p^{2} (k)<B \cdot (k+2)^{2} \cdot \ln ^{2} (k+2) {\rm \;} \text{ for all } k \in {\rm \mathbb{N}},\]

\noindent where $p^{2} (k)$ is a $T$-matrix leading element of $k$-row.

PROOF. It follows from Chebyshev inequalities \eqref{(26)} that there are 2 positive constants $c$ and $C$ such that
\[c\cdot n\cdot \ln n<p_{n} <C\cdot n\cdot \ln n {\rm \;}\text{ for all } n\ge 2,\]
\noindent where $p_{n} $ is the $n$-th prime number. In turn, the constants $c=\frac{1}{A} ,{\rm \; \; }C=\frac{2}{a} $ (see [4]) can be obtained from the proof of this proposition. Thereafter,
\begin{spacing}{-0.5}\end{spacing}
\[c\cdot n\cdot \ln n<p_{n} <C\cdot n\cdot \ln n,{\rm \; \; }n\ge 2{\rm \; \; }\Rightarrow {\rm \; \; }c\cdot n\cdot \ln n<p_{n} <C\cdot n\cdot \ln n,{\rm \; \; }n\ge 3{\rm \; \; }\Leftrightarrow\]
\begin{spacing}{-0.5}\end{spacing}
\[ \Leftrightarrow {\rm \; \; }c^{2} \cdot n^{2} \cdot \ln ^{2} n<p_{n}^{2} <C^{2} \cdot n^{2} \cdot \ln ^{2} n,{\rm \; \; }n\ge 3{\rm \; \; }\mathop{\Leftrightarrow }\limits^{{\rm \eqref{(2)}}} \] 
\[\Leftrightarrow {\rm \; \; }c^{2} \cdot (k+2)^{2} \cdot \ln ^{2} (k+2)<p^{2} (k)<C^{2} \cdot (k+2)^{2} \cdot \ln ^{2} (k+2),{\rm \; \; }k\in {\rm \mathbb{N}}  {\rm \; \; } \Leftrightarrow \] 
\[\Leftrightarrow {\rm \; \; }b\cdot (k+2)^{2} \cdot \ln ^{2} (k+2)<p^{2} (k)<B\cdot (k+2)^{2} \cdot \ln ^{2} (k+2),{\rm \;}k\in {\rm \mathbb{N}}, \text{where }  b\equiv c^{2} ,{\rm \; }B\equiv C^{2}.\]

Corollary 4.21 is proved.

PROPOSITION 4.22. 
\begin{spacing}{0.7}\end{spacing}\[1){\rm \; \;}2 \cdot \frac{\sqrt{x} }{\ln x}-2<\pi _{{\rm M} _{T} } (x)<2.51012 \cdot \frac{\sqrt{x} }{\ln x}-2,{\rm \; }x\ge 289.\]
\begin{spacing}{0.8}\end{spacing}
\[2){\rm \; \;} 2 \cdot \frac{\sqrt{x} }{\ln x}-2<\pi _{{\rm M} _{T} }(x) \le 2.2112...\cdot \frac{\sqrt{x} }{\ln x}-2 {\rm \;}\text{for sufficiently large $x$.}\]
\begin{spacing}{1.7}\end{spacing}

PROOF. 
$\pi (x)$ satisfies the inequalities (see [9]):
\begin{spacing}{0.8}\end{spacing}
\[0.9212...\cdot \frac{x}{\ln x} \le \pi (x)\le 1.1056...\cdot \frac{x}{\ln x} {\rm \;}\text{ for sufficiently large } x.\]
\begin{spacing}{1.5}\end{spacing}
Using the proof of Proposition 4.20, we make sure that with $a=0.9212...,{\rm \; }A=1.1056...$:
\[2 \cdot 0.9212... \cdot \frac{\sqrt{x} }{\ln x}-2 \le \pi _{{\rm M} _{T} } (x) \le 2 \cdot 1.1056...\cdot \frac{\sqrt{x} }{\ln x}-2 {\rm \; \;}\Leftrightarrow {\rm \; \;}\] 
\begin{equation} \label{(27)}
\Leftrightarrow {\rm \; \;}1.8424... \cdot \frac{\sqrt{x} }{\ln x}-2 \le \pi _{{\rm M} _{T} } (x) \le 2.2112...\cdot \frac{\sqrt{x} }{\ln x}-2
\end{equation} 
\begin{spacing}{1.5}\end{spacing}
\noindent for sufficiently large $x$. 

In turn, there is an inequality (see [10]):
\begin{spacing}{0.7}\end{spacing}
\[\frac{x}{\ln x} <\pi (x)<1.25506\cdot \frac{x}{\ln x} ,{\rm \; }x\ge 17.\]
\begin{spacing}{1.4}\end{spacing}

Use the proof of Proposition 4.20 again. As a result,
\begin{spacing}{0.8}\end{spacing}
\begin{equation} \label{(28)} 
2 \cdot \frac{\sqrt{x} }{\ln x}-2<\pi _{{\rm M} _{T} } (x)<2.51012 \cdot \frac{\sqrt{x} }{\ln x}-2,{\rm \; }x\ge 289.       
\end{equation}
\begin{spacing}{1.4}\end{spacing}

Then for sufficiently large $x$:
\begin{spacing}{0.8}\end{spacing}
\[2 \cdot \frac{\sqrt{x} }{\ln x}-2{\rm \;}\mathop{<}\limits^{\eqref{(28)}} {\rm \; }\pi _{{\rm M} _{T} } (x)\mathop{\le }\limits^{\eqref{(27)}}2.2112...\cdot \frac{\sqrt{x} }{\ln x}-2.\]
\begin{spacing}{1.4}\end{spacing}

Proposition 4.22 is proved.

THEOREM 4.23 (the asymptotic law of distribution of prime numbers).
\begin{spacing}{0.5}\end{spacing}\begin{equation} \label{(29)} 
\pi (x)\sim \frac{x}{\ln x} {\rm \; \; }(x\to \infty ),  
\end{equation} \begin{spacing}{1.4}\end{spacing}
\noindent where symbol $\sim$ stands for an asymptotic equality (an asymptotic equivalence):
\begin{spacing}{0.8}\end{spacing}\[f_{1} (x)\sim f_{2} (x){\rm \; \; }(x\to \infty ){\rm \; \; }\Leftrightarrow {\rm \; \; }\mathop{\lim }\limits_{x\to \infty } \frac{f_{1} (x)}{f_{2} (x)} =1 {\rm \; \; }\text{(see, for example, [7]).} \]
\begin{spacing}{1.5}\end{spacing}
PROPOSITION 4.24. 
\begin{spacing}{0.5}\end{spacing}
\begin{equation} \label{(30)}
\pi _{{\rm M} _{T} } \left(x\right)\sim \frac{2\sqrt{x} }{\ln x} {\rm \; \; }(x\to \infty ).  
\end{equation} 
\begin{spacing}{1.4}\end{spacing}

PROOF. By Theorem 4.23,  $\pi (x)\sim \frac{x}{\ln x} {\rm \; \; }(x\to \infty )$. That means that $\mathop{\lim }\limits_{x\to \infty } \frac{\pi (x)}{\left(\frac{x}{\ln x} \right)} =1$.  
\[\mathop{\lim }\limits_{x\to \infty } \frac{\pi _{{\rm M} _{T} } (x)}{\left(\frac{2\sqrt{x} }{\ln x} \right)} {\rm \; \; }\mathop{=}\limits^{\eqref{(25)}} {\rm \; }\mathop{\lim }\limits_{x\to \infty } \frac{\pi (\sqrt{x} )-2}{\left(\frac{\sqrt{x} }{\ln \sqrt{x} } \right)} =\left[{\rm \; }y=\sqrt{x} ,{\rm \; }y\mathop{\to }\limits_{x\to \infty } \infty \right]=\mathop{\lim }\limits_{y\to \infty } \frac{\pi (y)-2}{\left(\frac{y}{\ln y} \right)} =\] 
\[=\mathop{\lim }\limits_{y\to \infty } \frac{\pi (y)}{\left(\frac{y}{\ln y} \right)} -\mathop{\lim }\limits_{y\to \infty } \frac{2}{\left(\frac{y}{\ln y} \right)} =1-2\cdot \mathop{\lim }\limits_{y\to \infty } \frac{1}{\left(\frac{y}{\ln y} \right)} . \] 

Note that $\mathop{\lim }\limits_{y\to \infty } \frac{y}{\ln y} =\left[\frac{\infty }{\infty } \right]$. In this case, we use L'Hopital's rule to find the limit $\mathop{\lim }\limits_{y\to \infty } \frac{y}{\ln y}$.

\begin{spacing}{1.2}\end{spacing}\[\mathop{\lim }\limits_{y\to \infty } \frac{y}{\ln y} =\mathop{\lim }\limits_{y\to \infty } \frac{y'}{(\ln y)'} =\mathop{\lim }\limits_{y\to \infty } \frac{1}{\left(\frac{1}{y} \right)} =\infty {\rm \; \; }\Rightarrow {\rm \; \; }\mathop{\lim }\limits_{y\to \infty } \frac{1}{\left(\frac{y}{\ln y} \right)} ={\rm 0\; \; }\Rightarrow \]
\[\Rightarrow {\rm \; \; }\mathop{\lim }\limits_{x\to \infty } \frac{\pi _{{\rm M} _{T} } (x)}{\left(\frac{2\sqrt{x} }{\ln x} \right)} =1{\rm \; \; }\Leftrightarrow {\rm \; \; }\pi _{{\rm M} _{T} } (x)\sim \frac{2\sqrt{x} }{\ln x} {\rm \; \; }(x\to \infty ).\] 

Proposition 4.24 is proved.

COROLLARY 4.25. 
\begin{spacing}{0.5}\end{spacing}
\[\pi (x)\sim \frac{\sqrt{x} }{2} \cdot \pi _{{\rm M} _{T} } (x){\rm \; \; }(x\to \infty ).\]

PROOF. The asymptotic equality $\sim $ is an equivalence relation. Then, by the symmetric property of relation $\sim $,
\begin{spacing}{0.5}\end{spacing}\[\pi (x){\rm \; }\mathop{\sim }\limits^{(29)} {\rm \; }\frac{x}{\ln x} {\rm \; \; }(x\to \infty ){\rm \; \; }\Rightarrow {\rm \; \; }\frac{x}{\ln x} \sim {\rm \; }\pi (x){\rm \; \; }(x\to \infty ).\] 
\begin{spacing}{0.1}\end{spacing}
\[\frac{\sqrt{x} }{2} \cdot \pi _{{\rm M} _{T} } (x){\rm \; \; }\mathop{\sim }\limits^{(30)} {\rm \; \; }\frac{\sqrt{x} }{2} \cdot \left(\frac{2\sqrt{x} }{\ln x} \right){\rm \; \; }(x\to \infty ){\rm \; \; }\Leftrightarrow {\rm \; \; }\frac{\sqrt{x} }{2} \cdot \pi _{{\rm M} _{T} } (x)\sim \frac{x}{\ln x} {\rm \; \; }(x\to \infty ).\] 

By the transitive property of relation $\sim$,
\[ {\rm \; \; }\frac{\sqrt{x} }{2} \cdot \pi _{{\rm M} _{T} } (x)\sim \frac{x}{\ln x} {\rm \; \; }(x\to \infty ){\rm \; \;}\wedge {\rm \; \;}\frac{x}{\ln x} \sim {\rm \; }\pi (x){\rm \; \; }(x\to \infty ){\rm \;}\Rightarrow\]
\[\Rightarrow{\rm \;}\frac{\sqrt{x} }{2} \cdot \pi _{{\rm M} _{T} } (x)\sim {\rm \; }\pi (x){\rm \; \; }(x\to \infty ).\]
\begin{spacing}{1.5}\end{spacing}

Use the symmetric property of relation $\sim$ again, 
\begin{spacing}{0.7}\end{spacing}
\[\frac{\sqrt{x} }{2} \cdot \pi _{{\rm M} _{T} } (x)\sim {\rm \; }\pi (x){\rm \; \; }(x\to \infty ){\rm \; \; }\Rightarrow {\rm \; \; }\pi (x)\sim \frac{\sqrt{x} }{2} \cdot \pi _{{\rm M} _{T} } (x){\rm \; \; }(x\to \infty ).\] 
\begin{spacing}{1.4}\end{spacing}

Corollary 4.25 is proved. 

\newpage

\begin{center}
\textbf{\large 5. \textit{T}-matrix - based algorithm for finding all the prime numbers less than or equal to any given natural number \textit{n\ge 5}}
\end{center}
\begin{spacing}{0.8}\end{spacing}

The use of $T$-matrix lead to algorithm №1, similar to the sieve of Eratosthenes (see [11]) in asymptotic time complexity, for finding all the prime numbers less than or equal to a given natural number $n\ge 5$. We will present the description of this algorithm in a proof of next theorem on the basis of $T$-matrix. 

Further, let $a\% b$ denote the remainder after dividing $a\in {\rm \mathbb{N}}$ by $b\in {\rm \mathbb{N}}$.

THEOREM 5.1.
\begin{spacing}{0.3}\end{spacing}\[\left(\forall b\not \in {\rm \mathbb{P}} \right)\left(\left(\exists h\in {\rm \mathbb{N}} \right)\left(b=6h+1{\rm \;}\vee {\rm \;}b=6h-1 \right){\rm \;}\Rightarrow {\rm \;} b \in \widetilde{T} \right).\]
\begin{spacing}{1.3}\end{spacing}

PROOF. \textbf{Description of algorithm №1. Input:  } $n\in {\rm \mathbb{N}}$: $n\ge 5$. Initially, let's say that 2 and 3 are prime numbers, since they could not be represented as $6r\pm 1$, $r\in {\rm \mathbb{N}} $. 

\textbf{Step 1. }Write consecutively all numbers $6r-1;{\rm \; }6r+1$ less than or equal to given number $n$, where $r\in {\rm \mathbb{N}} $ is a number of pair. The received list of numbers denote by $pr$. Initially, let a variable $p$ equals to $p(1)$, that is the third prime number.

\textbf{Step 2. }Verify the conditions:

1) If $p^{2} >x$, then the algorithm halts. \textbf{Output: } list of all the prime numbers less than or equal to $n$.

2) If $p\% 6=5$, then cross out in the list $pr$ the numbers, starting at the leading element $p^{2}$:
\begin{spacing}{0.3}\end{spacing}\begin{equation} \label{(31)} 
p^{2} , p^{2} +2p, p^{2} +2p+4p, p^{2} +2p+4p+2p, p^{2} +2p+4p+2p+4p, \dots \le n.  
\end{equation} \begin{spacing}{0.8}\end{spacing}

All crossed out numbers \eqref{(31)} are increased by one in the list $pr$. Proceed to Step 3.

3) If $p\% 6=1$, then cross out in the list $pr$ the numbers, starting at the leading element $p^{2}$:
\begin{spacing}{0.3}\end{spacing}\begin{equation} \label{(32)} 
p^{2} , p^{2} +4p, p^{2} +4p+2p, p^{2} +4p+2p+4p, p^{2} +4p+2p+4p+2p, \dots \le n.  
\end{equation}
\begin{spacing}{0.8}\end{spacing}

All crossed out numbers \eqref{(32)} are increased by one in the list $pr$. Proceed to step 3.

 \textbf{Step 3. } Find the first uncrossed out number greater than $p$ in the list $pr$, and let $p$ now equal this new number. Proceed to Step 2.
 
 As a result, we get that the composite numbers, obtained in the $k$-th deletion of Step 2, are some elements of sequence $T_{k} {\rm \;}(k\in \mathbb{N})$. Thus, all elements of sequences $T_{k}$ for any $k\in {\rm \mathbb{N}}$ will be crossed out as $n\to +\infty $. That means that all composite numbers of the form $6h\pm 1$ fully exhaust all elements of $T$-matrix. 
 
Theorem 5.1 is proved.
\begin{spacing}{1.5}\end{spacing}
The composite numbers are shown in red, everyone else are prime numbers in Table №2. Each cell of Table №2 stores the number pair $(6h-1;{\rm \; }6h+1)$. The number $h\in {\rm \mathbb{N}}$ of the pair $(6h-1;{\rm \; }6h+1)$ is specified in the top-left corner of the cell. Each composite number has index $k\in {\rm \mathbb{N}}$ next to it. This index determines at what stage $k$ of Step 2 of algorithm №1 the composite number was crossed out. The composite numbers, encountered only 1 time in $T$-matrix, are circled. They are located in the first row of $T$-matrix. 
\[\textbf{Table №2}.\text{ The obtaining all the prime numbers } (2, 3 \in \mathbb{P}) \text{ less than or equal to 508}\]

\begin{center}\includegraphics*[scale=0.6]{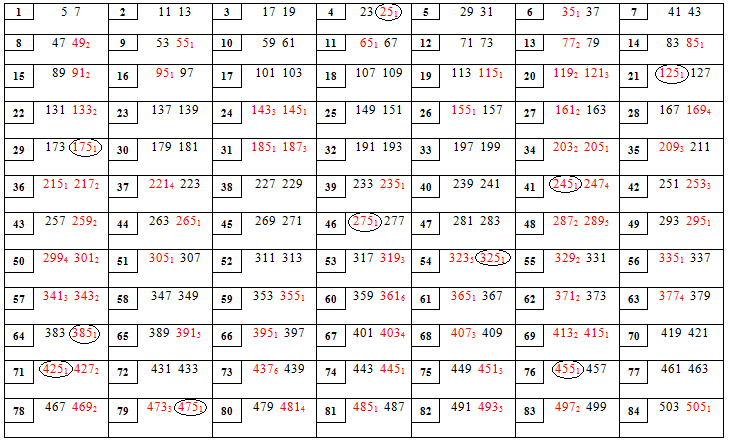} \end{center}

Further, program implementation of algorithm №1 in Java (see «Source code №1») and output at $n=5000$ will be presented.

\begin{spacing}{-0.5}\end{spacing}
\begin{center}\textbf{Source code №1.} Implementation of algorithm №1 on the basis of $T$-matrix to obtain all the prime numbers less than or equal to given natural number $n\ge 5$
\end{center}

\begin{center}\begin{tabular}{|p{6in}|} \hline 
Sieve.java \\ \hline 
\small \textcolor[RGB]{127,0,85} {\textbf{public class}} Sieve

\{
\begin{spacing}{0.9}\end{spacing} 
{\rm \; \; \; \;}\textcolor {SeaGreen}{// function defining the elements of sequence }

{\rm \; \; \;} \textcolor {SeaGreen}{// 5, 7, 11, 13, 17, 19, 23, 25, ... , 6k-1, 6k+1,...}

 {\rm \; \; \;} \textcolor[RGB]{127,0,85} {\textbf{public static int}} f(\textcolor[RGB]{127,0,85}{\textbf{int}} i) \{\textcolor[RGB]{127,0,85}{\textbf{return}} (\textcolor[RGB]{127,0,85}{\textbf{int}})(3*i+1.5-0.5*Math.\textit{pow}(-1,i));\}

 {\rm \; \; \;} \textcolor {SeaGreen}{// function defining the elements of sequence}
 
 {\rm \; \; \;} \textcolor {SeaGreen}{// p\^{}2, p\^{}2+2p, p\^{}2+6p, p\^{}2+8p, p\^{}2+12p, ...}
 
 {\rm \; \; \;} \textcolor[RGB]{127,0,85}{\textbf{public static int}} J1(\textcolor[RGB]{127,0,85}{\textbf{int}} p, \textcolor[RGB]{127,0,85} {\textbf{int}} j) \{\textcolor[RGB]{127,0,85}{\textbf{return}} (\textcolor[RGB]{127,0,85}{\textbf{int}})(p*p+(\textit{f}(j)-5)*p);\}
 
 {\rm \; \; \;} \textcolor {SeaGreen}{// function defining the elements of sequence}

{\rm \; \; \;}\textcolor {SeaGreen}{// p\^{}2, p\^{}2+4p, p\^{}2+6p, p\^{}2+10p, p\^{}2+12p, ...}

 {\rm \; \; \;} \textcolor[RGB]{127,0,85}{\textbf{public static int}} J2(\textcolor[RGB]{127,0,85}{\textbf{int}} p, \textcolor[RGB]{127,0,85}{\textbf{int}} j) \{\textcolor[RGB]{127,0,85}{\textbf{return}} (\textcolor[RGB]{127,0,85}{\textbf{int}})(p*p+(\textit{f}(j)+Math.\textit{pow}(-1, j)-4)*p);\}
 
{\rm \; \; \;} \textcolor {SeaGreen}{// function defining the index of a given element in a list by using binary search technique}

{\rm \; \; \;} \textcolor[RGB]{127,0,85} {\textbf{public static int}} binarySearch(\textcolor[RGB]{127,0,85}{\textbf{int}}[] a, \textcolor[RGB]{127,0,85}{\textbf{int}} key)  
 
  {\rm \; \; \;} \{
 
 {\rm \; \; \; \; \; \;}\textcolor[RGB]{127,0,85}{\textbf{int}} min\_index = 0, max\_index = a.\textcolor{Blue}{length} - 1, mid\_index; 
\\ \hline 
\end{tabular}\end{center}

\begin{center}\begin{tabular}{|p{6in}|}
\hline 
{\rm \; \; \; \; \; \;}\textcolor[RGB]{127,0,85}{\textbf{while}} (min\_index $<$= max\_index) 

{\rm \; \; \; \; \; \;}\{          
 
{\rm \; \; \; \; \; \; \; \; \;}\textcolor {SeaGreen}{// index of the middle element of a list}
 
 {\rm \; \; \; \; \; \; \; \; \;}mid\_index = (min\_index + max\_index) $>$$>$$>$ 1;
 
{\rm \; \; \; \; \; \; \; \; \;}\textcolor[RGB]{127,0,85}{\textbf{if}} (a[mid\_index] $<$ key)

{\rm \; \; \; \; \; \; \; \; \; \; \; \;}min\_index = mid\_index + 1;

{\rm \; \; \; \; \; \; \; \; \;}\textcolor[RGB]{127,0,85}{\textbf{else if}} (a[mid\_index] $>$ key)

{\rm \; \; \; \; \; \; \; \; \; \; \; \;} max\_index = mid\_index - 1; 

{\rm \; \; \; \; \; \; \; \; \;}\textcolor[RGB]{127,0,85}{\textbf{else return}} mid\_index; \textcolor {SeaGreen}{//  index of a given element is found in a list} 

{\rm \; \; \; \; \; \;}\}

{\rm \; \; \; \; \; \;}\textcolor[RGB]{127,0,85}{\textbf{return}} 0; \textcolor {SeaGreen}{// index of a given element is not found in a list}

{\rm \; \; \; }\}

{\rm \; \; \; } \textcolor[RGB]{127,0,85}{\textbf{public static void}} main(String[] args)

{\rm \; \; \; }\{

{\rm \; \; \;\; \; \; }\textcolor {SeaGreen}{//  input}

{\rm \; \; \;\; \; \; }\textcolor[RGB]{127,0,85}{\textbf{int}} n=5000,

\small {\rm \; \; \; \; \; \;}\textcolor {SeaGreen}{// counters of loops}

{\rm \; \; \;\; \; \; }i=1, j,

{\rm \; \; \;\; \; \; }\textcolor {SeaGreen}{// uncrossed out prime number}

{\rm \; \; \;\; \; \; }p; 

{\rm \; \; \; \; \; \;}\textcolor[RGB]{127,0,85}{\textbf{if}} (n$<$=1) \textcolor[RGB]{127,0,85}{\textbf{return}};

{\rm \; \; \; \; \; \;}\textcolor[RGB]{127,0,85}{\textbf{if}} (n==2) \{System.\textcolor{Blue}{\textit{out}}.print(\textcolor{Blue}{''2''}); \textcolor[RGB]{127,0,85}{\textbf{return}};\}
   
{\rm \; \; \; \; \; \;}\textcolor[RGB]{127,0,85}{\textbf{if}} (n$<$=4) \{System.\textcolor{Blue}{\textit{out}}.print(\textcolor{Blue}{''2, 3''}); \textcolor[RGB]{127,0,85}{\textbf{return}};\}

{\rm \; \; \; \; \; \; }\textcolor {SeaGreen}{// \textit{pr} is a list to store 0 and the numbers of the form 6k-1,6k+1} 

{\rm \; \; \; \; \; \; }\textcolor[RGB]{127,0,85}{\textbf{int}} [] pr=\textcolor[RGB]{127,0,85}{\textbf{new int}}[(n+2)/3-(n\%6)/4+(n\%6)/5];

{\rm \; \; \; \; \; \; }\textcolor {SeaGreen}{// fill the list \textit{pr}}

{\rm \; \; \; \; \; \; }\textcolor[RGB]{127,0,85}{\textbf{while}} (\textit{f}(i)$<$=n)     

{\rm \; \; \;\; \; \; }\{

{\rm \; \; \; \; \; \; \; \; \;}pr[i]=\textit{f}(i);
 
{\rm \; \; \; \; \; \; \; \; \;}i++;

{\rm \; \; \; \; \; \; }\}   
    
{\rm \; \; \; \; \; \; }\textcolor[RGB]{127,0,85}{\textbf{for}} (i=1; pr[i]*pr[i]$<$=n; i++) 

{\rm \; \; \; \; \; \; \; \; \;  }\textcolor {SeaGreen}{// condition: a number is not crossed out (prime) in the list \textit{pr}}

{\rm \; \; \; \; \; \; \; \; \;  }\textcolor[RGB]{127,0,85}{\textbf{if}} ((pr[i]\%2)!=0)

{\rm \; \; \; \; \; \; \; \; \; }\{

{\rm \; \; \; \; \; \; \; \; \; \; \; \;  }p=pr[i]; 

{\rm \; \; \; \; \; \; \; \; \; \; \; \;  }\textcolor[RGB]{127,0,85}{\textbf{if}} (p\%6==5)

{\rm \; \; \; \; \; \; \; \; \; \; \; \; \; \; \;  }\textcolor[RGB]{127,0,85}{\textbf{for}} (j=1; \textit{J1}(p, j)$<$=n; j++)
\\ \hline 
\end{tabular}\end{center}

\begin{center}\begin{tabular}{|p{6in}|}
\hline 
{\rm \; \; \; \; \; \; \; \; \; \; \; \;\; \; \; \; \; \; }\textcolor {SeaGreen}{// deletion of the composite numbers}

{\rm \; \; \; \; \; \; \; \; \; \; \; \;\; \; \; \; \; \; }\textcolor {SeaGreen}{// p\^{}2, p\^{}2+2p, p\^{}2+6p, p\^{}2+8p, p\^{}2+12p, ...$<$=n}

{\rm \; \; \; \; \; \; \; \; \; \; \; \;\; \; \; \; \; \; }pr[\textit{binarySearch}(pr, \textit{J1}(p, j))]++;

{\rm \; \; \; \; \; \; \; \; \;  \; \; \; }\textcolor[RGB]{127,0,85}{\textbf{if}} (p\%6==1)

{\rm \; \; \; \; \; \; \; \; \; \; \; \; \; \; \; }\textcolor[RGB]{127,0,85}{\textbf{for}} (j=1; \textit{J2}(p, j)$<$=n; j++)

{\rm \; \; \; \; \; \; \; \; \; \; \; \;\; \; \; \; \; \; }\textcolor {SeaGreen}{//  deletion of the composite numbers}

{\rm \; \; \; \; \; \; \; \; \; \; \; \;\; \; \; \; \; \; }\textcolor {SeaGreen}{// p\^{}2, p\^{}2+4p, p\^{}2+6p, p\^{}2+10p, p\^{}2+12p, ...$<$=n} 

{\rm \; \; \; \; \; \; \; \; \; \; \; \;\; \; \; \; \; \; }pr[\textit{binarySearch}(pr, \textit{J2}(p, j))]++;
 
{\rm \; \; \; \; \; \; \; \; \; }\}
 
{\rm \; \; \; \; \; \; }\textcolor {SeaGreen}{// print the prime numbers 2 and 3}

{\rm \; \; \; \; \; \; }System.\textcolor {Blue}{\textit out}.print(\textcolor {Blue}{''2, 3''});

{\rm \; \; \; \; \; \; }\textcolor[RGB]{127,0,85}{\textbf{for}} (i=1; i$<$pr.\textcolor {Blue}{length};  i++)

{\rm \; \; \; \; \; \; \; \; \; }\textcolor[RGB]{127,0,85}{\textbf{if}} ((pr[i]\%2)!=0)

{\rm \; \; \; \; \; \; \; \; \; \; \; \;}\textcolor {SeaGreen}{// print the prime numbers between 5 and n}

{\rm \; \; \; \; \; \; \; \; \; \; \; \;}System.\textcolor{Blue}{\textit{out}}.print(\textcolor {Blue}{'',''}+pr[i]);

{\rm \; \; \; }\}

\}

\small // Output:\newline \scriptsize {2,3,5,7,11,13,17,19,23,29,31,37,41,43,47,53,59,61,67,71,73,79,83,89,97,101,103,107,109,113,127,131,137,139,149,151,157,
163,167,173,179,181,191,193,197,199,211,223,227,229,233,239,241,251,257,263,269,271,277,281,283,293,307,311,313,317,
331,337,347,349,353,359,367,373,379,383,389,397,401,409,419,421,431,433,439,443,449,457,461,463,467,479,487,491,499,
503,509,521,523,541,547,557,563,569,571,577,587,593,599,601,607,613,617,619,631,641,643,647,653,659,661,673,677,683,
691,701,709,719,727,733,739,743,751,757,761,769,773,787,797,809,811,821,823,827,829,839,853,857,859,863,877,881,883,
887,907,911,919,929,937,941,947,953,967,971,977,983,991,997,1009,1013,1019,1021,1031,1033,1039,1049,1051,1061,1063,
1069,1087,1091,1093,1097,1103,1109,1117,1123,1129,1151,1153,1163,1171,1181,1187,1193,1201,1213,1217,1223,1229,
1231,1237,1249,1259,1277,1279,1283,1289,1291,1297,1301,1303,1307,1319,1321,1327,1361,1367,1373,1381,1399,1409,
1423,1427,1429,1433,1439,1447,1451,1453,1459,1471,1481,1483,1487,1489,1493,1499,1511,1523,1531,1543,1549,1553,
1559,1567,1571,1579,1583,1597,1601,1607,1609,1613,1619,1621,1627,1637,1657,1663,1667,1669,1693,1697,1699,1709,
1721,1723,1733,1741,1747,1753,1759,1777,1783,1787,1789,1801,1811,1823,1831,1847,1861,1867,1871,1873,1877,1879,
1889,1901,1907,1913,1931,1933,1949,1951,1973,1979,1987,1993,1997,1999,2003,2011,2017,2027,2029,2039,2053,2063,
2069,2081,2083,2087,2089,2099,2111,2113,2129,2131,2137,2141,2143,2153,2161,2179,2203,2207,2213,2221,2237,2239,
2243,2251,2267,2269,2273,2281,2287,2293,2297,2309,2311,2333,2339,2341,2347,2351,2357,2371,2377,2381,2383,2389,
2393,2399,2411,2417,2423,2437,2441,2447,2459,2467,2473,2477,2503,2521,2531,2539,2543,2549,2551,2557,2579,2591,
2593,2609,2617,2621,2633,2647,2657,2659,2663,2671,2677,2683,2687,2689,2693,2699,2707,2711,2713,2719,2729,2731,
2741,2749,2753,2767,2777,2789,2791,2797,2801,2803,2819,2833,2837,2843,2851,2857,2861,2879,2887,2897,2903,2909,
2917,2927,2939,2953,2957,2963,2969,2971,2999,3001,3011,3019,3023,3037,3041,3049,3061,3067,3079,3083,3089,3109,
3119,3121,3137,3163,3167,3169,3181,3187,3191,3203,3209,3217,3221,3229,3251,3253,3257,3259,3271,3299,3301,3307,
3313,3319,3323,3329,3331,3343,3347,3359,3361,3371,3373,3389,3391,3407,3413,3433,3449,3457,3461,3463,3467,3469,
3491,3499,3511,3517,3527,3529,3533,3539,3541,3547,3557,3559,3571,3581,3583,3593,3607,3613,3617,3623,3631,3637,
3643,3659,3671,3673,3677,3691,3697,3701,3709,3719,3727,3733,3739,3761,3767,3769,3779,3793,3797,3803,3821,3823,
3833,3847,3851,3853,3863,3877,3881,3889,3907,3911,3917,3919,3923,3929,3931,3943,3947,3967,3989,4001,4003,4007,
4013,4019,4021,4027,4049,4051,4057,4073,4079,4091,4093,4099,4111,4127,4129,4133,4139,4153,4157,4159,4177,4201,
4211,4217,4219,4229,4231,4241,4243,4253,4259,4261,4271,4273,4283,4289,4297,4327,4337,4339,4349,4357,4363,4373,
4391,4397,4409,4421,4423,4441,4447,4451,4457,4463,4481,4483,4493,4507,4513,4517,4519,4523,4547,4549,4561,4567,
4583,4591,4597,4603,4621,4637,4639,4643,4649,4651,4657,4663,4673,4679,4691,4703,4721,4723,4729,4733,4751,4759,
4783,4787,4789,4793,4799,4801,4813,4817,4831,4861,4871,4877,4889,4903,4909,4919,4931,4933,4937,4943,4951,4957,
4967,4969,4973,4987,4993,4999} \\ \hline 
\end{tabular}\end{center}

\textbf{}

THEOREM 5.2 (Mertens formula). 
\begin{equation} \label{(33)} 
\sum _{p\in {\rm \mathbb{P}} :{\rm \; }p\le x}\frac{1}{p}=M+\ln (\ln x)+O\left(\frac{1}{\ln x} \right), x>2,
\end{equation}
\noindent where $M$ is the Mertens constant independent of $x$: $M\approx 0.2614972128476427$ (see [6]).
 
 For practical computation of constant $M$ we can use relationship (see [7])
\[M=\gamma +\sum _{p\in {\rm \mathbb{P}} }\left(\ln \left(1-\frac{1}{p} \right)+\frac{1}{p} \right) ,\] 

\noindent where $\gamma $ is the Euler-Mascheroni constant: $\gamma \approx 0.577215664901533$.

PROPOSITION 5.3. The asymptotic time complexity of the algorithm №1 is $O(n\ln (\ln n))$.

RROOF. Let $pr$ is a list comprising numbers of the form $6h\pm 1$, less than or equal to a given natural number $n\ge 5$ (as in algorithm №1, Step 1). The numbering of the pairs $(6h-1;{\rm \; }6h+1)$ starts with $h=1$. For convenience, let's say that the null element of the list $pr$ is 0. A number $\tilde{m}_{k} (n)\in {\rm \mathbb{N}}$ of composites, less than or equal to given number $n$, can be estimated based on equality \eqref{(3)} in the $k$-th deletion ($k\in {\rm \mathbb{N}} $) of Step 2 (further, in Step 2,$k$):
\[p(k)\cdot \left(3\tilde{m}_{k} (n)+\frac{3-(-1)^{\tilde{m}_{k} (n)} }{2} \right)\le n{\rm \; \; \; }\Leftrightarrow {\rm \; \; \; }3\tilde{m}_{k} (n)+\frac{3-(-1)^{\tilde{m}_{k} (n)} }{2} \le \frac{n}{p(k)} {\rm \; \; \; }\Leftrightarrow \] 
\[\Leftrightarrow {\rm \; \; \; }\tilde{m}_{k} (n)+\frac{3-(-1)^{\tilde{m}_{k} (n)} }{6} \le \frac{n}{3\cdot p(k)} {\rm \; \; \; }\Rightarrow {\rm \; \; \; }\tilde{m}_{k} (n)\le \frac{n}{3\cdot p(k)} -\frac{3-(-1)^{\tilde{m}_{k} (n)} }{6} =\]
\begin{spacing}{0.5}\end{spacing}
\[=\frac{n}{3\cdot p(k)} -\frac{1}{2} +\frac{(-1)^{\tilde{m}_{k} (n)} }{6} \le \frac{n}{3\cdot p(k)} -\frac{1}{2} +\frac{1}{6} =\frac{n}{3\cdot p(k)} -\frac{1}{3} <\frac{n}{3\cdot p(k)} {\rm \; \; \; }\Rightarrow \] 
\begin{equation} \label{(34)} 
\Rightarrow{\rm \; \; }\tilde{m}_{k} (n)<\frac{n}{3\cdot p(k)} .     
\end{equation} \begin{spacing}{1.2}\end{spacing}

In Step 2,$k$, the composite numbers are crossed out starting with leading element $p^{2} (k)$ in $k$-row of $T$-matrix, since all composite numbers less than $p^{2} (k)$ are already crossed out in Steps 2,$r$, where $r=\overline{1;k-1}$. Therefore, in each Step 2,$k$, a number $\Delta _{k} \in {\rm\mathbb{N}} _{0}$ of uncrossed out composites less than $p^{2} (k)$ is defined by the formula:
\begin{spacing}{0.5}\end{spacing}\begin{equation} \label{(35)} 
\Delta _{k} ={\rm \; }\# _{k} (p^{2} (k))-1,  
\end{equation}
\noindent where  $\# _{k} (p^{2} (k))$ is a number of leading element $p^{2} (k)$ in $k$-row of $T$-matrix. It is clear that $\Delta _{k} $ does not depend on the number $n$. In turn, the time complexity of binary search for defining the number $\# _{k} (p^{2} (k))$ (the index of element $p^{2} (k)$ in the list $pr$) is $O(\log _{2} C(n))$, where $C(n)$ is a space complexity of algorithm №1 for given number $n\in\mathbb{N}$: $n\ge 5$, due to dividing the list $pr$ in half (about binary search see [8]). In each Step 2,$k$, a number $m_{k} (n)$ of crossed out composites less than or equal to $n$: 
\begin{spacing}{0.3}\end{spacing}\begin{equation} \label{(36)} 
m_{k} (n)=\tilde{m}_{k} (n)-\Delta _{k} .                                      
\end{equation} 

 We will estimate the aggregate number $\sum _{k=1}^{\pi _{{\rm M} _{T} } (n)}m_{k} (n)$ of crossed out composites, less than or equal to $n$, as follows:
\begin{spacing}{0.8}\end{spacing}\[\sum _{k=1}^{\pi _{{\rm M} _{T} } (n)}m_{k} (n) {\rm \;} \mathop{=}\limits^{\eqref{(36)}} {\rm \;}\sum _{k=1}^{\pi _{{\rm M} _{T} } (n)}(\tilde{m}_{k} (n)-\Delta _{k} ) \le \sum _{k=1}^{\pi _{{\rm M} _{T} } (n)}\tilde{m}_{k} (n){\rm \; }\mathop{<}\limits^{\eqref{(34)}} {\rm \;} \sum _{k=1}^{\pi _{{\rm M} _{T} } (n)}\frac{n}{3\cdot p(k)}=\]
\begin{spacing}{0.5}\end{spacing}
\[=\frac{n}{3} \cdot \sum _{k=1}^{\pi _{{\rm M} _{T}} (n)}\frac{1}{p(k)} \le \frac{n}{3} \cdot \sum _{k=1}^{\pi _{{\rm M} _{T}} (n+4)}\frac{1}{p(k)} {\rm \; }\mathop{=}\limits^{\eqref{(25)}}  {\rm \; } \frac{n}{3} \cdot \sum _{k=1}^{\pi \left( \sqrt{n+4} \right)-2}\frac{1}{p(k)}{\rm \; }<\frac{n}{3} \cdot \sum _{k=1}^{\pi \left( \sqrt{n+4} \right)}\frac{1}{p(k)}=\]
\begin{spacing}{0.5}\end{spacing}
\[=\frac{n}{3} \cdot \left(\sum _{k=1}^{\pi \left( \sqrt{n+4} \right)}\frac{1}{p(k)}  +\frac{1}{2} -\frac{1}{2} +\frac{1}{3} -\frac{1}{3} \right)=\frac{n}{3} \cdot \left(\left(\sum _{k=1}^{\pi \left( \sqrt{n+4} \right)}\frac{1}{p(k)}  +\frac{1}{2} +\frac{1}{3} \right)-\frac{1}{2} -\frac{1}{3} \right)=\]

\[{\rm \; \;} =\frac{n}{3} \cdot \left(\sum _{p\in {\rm \mathbb{P}} :{\rm \; }p\le \sqrt{n+4} }\frac{1}{p} -\frac{5}{6}  \right)<\frac{n}{3} \cdot \sum _{p\in {\rm \mathbb{P}} :{\rm \; }p\le \sqrt{n+4}}\frac{1}{p} {\rm \; \;} \mathop{=}\limits^{\eqref{(33)}}\]
 \begin{spacing}{1.5}\end{spacing} 
\[= \frac{n}{3} \cdot \left(M+ \ln \left(\ln \sqrt{n+4} \right)+O\left(\frac{1}{\ln \sqrt{n+4}}\right)\right)=\] 
 \begin{spacing}{1.3}\end{spacing} 
\[=\frac{n}{3} \cdot M+ \frac{n}{3} \cdot \ln \left(\ln \sqrt{n+4} \right)+\frac{n}{3} \cdot O\left(\frac{1}{\ln \sqrt{n+4}}\right).\]
 \begin{spacing}{1.3}\end{spacing} 
\[\ln \sqrt{n+4} =\frac{1}{2} \cdot \ln (n+4)<\frac{1}{2} \cdot \ln \left(n^2 \right)=\ln n, n \ge 5 {\rm \; \;}\Rightarrow\]
 \begin{spacing}{0.9}\end{spacing} 
\begin{equation} \label{(37)} 
\Rightarrow {\rm \;}\sum _{k=1}^{\pi _{{\rm M} _{T} } (n)}m_{k} (n)=O(n \ln (\ln n)).     
\end{equation} \begin{spacing}{1.2}\end{spacing} 

 In Corollary 5.6 from Proposition 5.5, the asymptotic space complexity of the algorithm №1 is $O(n)$. Thus, we get for the time complexity $m(n)$ of algorithm №1:
\begin{spacing}{0.3}\end{spacing}\[m(n)=\sum _{k=1}^{\pi _{{\rm M} _{T} } (n)}\left(m_{k} (n)+O(\log _{2} C(n))\right) =\sum _{k=1}^{\pi _{{\rm M} _{T} } (n)}m_{k} (n) +\sum _{k=1}^{\pi _{{\rm M} _{T} } (n)}O(\log _{2} C(n)) {\rm \; }=\]
\[{\rm \; \; \; \; }=\sum _{k=1}^{\pi _{{\rm M} _{T} } (n)}m_{k} (n) +\sum _{k=1}^{\pi _{{\rm M} _{T} } (n)}O(\log _{2} n) \mathop{=}\limits^{\eqref{(37)}}O(n \ln (\ln n))+\sum _{k=1}^{\pi _{{\rm M} _{T} } (n)}O(\log _{2} n)=\]
\[=O(n \ln (\ln n))+O(\log _{2} n) \cdot \sum _{k=1}^{\pi _{{\rm M} _{T}} (n)} 1=O(n\ln (\ln n))+ \pi _{{\rm M} _{T}} (n) \cdot O(\log _{2} n).\]

Using the upper estimate from Proposition 4.20 for $\pi _{{\rm M} _{T}} (n+4)$,
\begin{spacing}{1.2}\end{spacing}
\[ \pi _{{\rm M} _{T}} (n) \le \pi _{{\rm M} _{T}} (n+4)=O \left(\frac{\sqrt{n+4}}{\ln (n+4)}\right).\]
\begin{spacing}{1.2}\end{spacing}
\[ \frac{\sqrt{n+4}}{\ln (n+4)}<\frac{\sqrt{2 \cdot n}}{\ln (n+4)}<\sqrt{2} \cdot \left(\frac{\sqrt{n}}{\ln n}\right), n \ge 5 {\rm \; \; }\Rightarrow\]
\[\Rightarrow {\rm \; \; }\pi _{{\rm M} _{T}} (n)=O\left( \frac{\sqrt{n}}{\ln n} \right){\rm \;}\Rightarrow{\rm \;}m(n)=O(n \ln (\ln n))+ O\left( \frac{\sqrt{n}}{\ln n} \right) \cdot O(\ln n)=\]
\begin{spacing}{1.3}\end{spacing}
\[{\rm \; \; \; \; }=O(n\ln (\ln n))+ O\left( \frac{\sqrt{n}}{\ln n} \cdot \ln n \right)=O(n \ln (\ln n))+ O(\sqrt{n}){\rm \; \; }\Rightarrow\]
\begin{spacing}{0.8}\end{spacing}
\[ \Rightarrow{\rm \; \; } m(n)=O(n \ln (\ln n)).\] 

 Proposition 5.3 is proved.

PROPOSITION 5.4. For any given natural number $n\ge 5$, in each Step 2,$k$ of algorithm №1:

 1) if $p(k)\% 6=5$, then in $k$-row of $T$-matrix, the crossed out elements have the form
\begin{spacing}{0.5}\end{spacing}\begin{equation} \label{(38)} 
p^{2} (k)+(f(j)-5)\cdot p(k),{\rm \; }j=\overline{1;m_{k} (n)}.     
\end{equation} 

2) if $p(k)\% 6=1$, then in $k$-row of $T$-matrix, the crossed out elements have the form
\begin{spacing}{0.5}\end{spacing}
\begin{equation} \label{(39)} 
p^{2} (k)+(f(j)+(-1)^{j} -4)\cdot p(k),{\rm \; }j=\overline{1;m_{k} (n)}.   
\end{equation} 
\noindent Here $m_{k} (n)$ is a number of crossed out composites less than or equal to $n$ in Step 2,$k$.

PROOF. 1) Let $p(k)\%6=5$. Then it follows from Step 2,$k$, para 2) of algorithm №1 that 
\begin{spacing}{-0.5}\end{spacing}
\[p^{2} (k),{\rm \; }p^{2} (k)+2 \cdot p(k),{\rm \; }p^{2} (k)+6 \cdot p(k),{\rm \; }p^{2} (k)+8 \cdot p(k),{\rm \; }p^{2} (k)+12 \cdot p(k),{\rm \; }...{\rm \; },{\rm \; }p^{2} (k)+a_{m_{k} (n)} \cdot p(k)\le n\]
\begin{spacing}{0.5}\end{spacing}
\noindent are composite numbers. As a result, we get a finite sequence of the form
\begin{spacing}{-0.3}\end{spacing}
\begin{equation} V\equiv {\rm \; (}a_{j} {\rm )}_{j=1}^{m_{k} (n)}: 0;2;6;8;12;...;a_{m_{k} (n)}, \text{where } a_{m_{k} (n)} \in {\rm \mathbb{N}}: p^{2} (k)+a_{m_{k} (n)} \cdot p(k)\le n.\end{equation}
\begin{spacing}{0.8}\end{spacing} By analogy with construction of sequence $T_{k} $ (see (1)), we find the formula for calculation of elements of sequence $V$:
\begin{spacing}{0.8}\end{spacing}\[a_{j} =2\cdot \left(\left\lfloor \frac{j}{2} \right\rfloor +2\cdot \left\lfloor \frac{j-1}{2} \right\rfloor \right),{\rm \; }j=\overline{1;m_{k} (n)}.\] 
\begin{spacing}{1.3}\end{spacing}
\[a_{j} {\rm \; }\mathop{=}\limits^{\eqref{(6)}} {\rm \; }2\cdot \left(\left\lfloor \frac{j}{2} \right\rfloor +2\cdot \left(j-1-\left\lfloor \frac{j}{2} \right\rfloor \right)\right)=2\cdot \left(\left\lfloor \frac{j}{2} \right\rfloor +2j-2-2\cdot \left\lfloor \frac{j}{2} \right\rfloor \right)=\]
\[=2\cdot \left(2j-\left\lfloor \frac{j}{2} \right\rfloor -2\right)=2\cdot \left(2j-\left(\frac{j}{2} -\frac{1-(-1)^{j} }{4} \right)-2\right)=2\cdot \left(2j-\frac{j}{2} +\frac{1-(-1)^{j} }{4} -2\right)=\] 
\begin{spacing}{-0.1}\end{spacing}
\[=3j+\frac{1-(-1)^{j} }{2}-4=\left(3j+\frac{3-(-1)^{j} }{2} \right)-1-4=f(j)-5,{\rm \; }j=\overline{1;m_{k} (n)}.\] 

Therefore, in $k$-row of $T$-matrix, the crossed out elements have form \eqref{(38)}.

2) Let $p(k)\% 6=1$. Then it follows from Step 2,$k$, para 3) of algorithm №1 that
\begin{spacing}{-0.5}\end{spacing}
\[p^{2} (k),{\rm \; }p^{2} (k)+4 \cdot p(k),{\rm \; }p^{2} (k)+6 \cdot p(k),{\rm \; }p^{2} (k)+10 \cdot p(k),{\rm \; }p^{2} (k)+12 \cdot p(k),...,p^{2} (k)+b_{m_{k} (n)} \cdot p(k)\le n\]
\begin{spacing}{0.8}\end{spacing}
\noindent are composite numbers. As a result, we get a finite sequence of the form
\begin{spacing}{-0.3}\end{spacing}
\begin{equation} W\equiv (b_{j} )_{j=1}^{m_{k} (n)}: 0;4;6;10;12;...;b_{m_{k} (n)}, \text{where } b_{m_{k} (n)} \in {\rm \mathbb{N}}: p^{2} (k)+b_{m_{k} (n)} \cdot p(k)\le n.\end{equation}
\begin{spacing}{0.8}\end{spacing}
By analogy with construction of sequence $T_{k}$ (see (1)), we find the formula for calculation of elements of sequence $W$:
\begin{spacing}{0.9}\end{spacing}
\[b_{j} =2\cdot \left(2\cdot \left\lfloor \frac{j}{2} \right\rfloor +\left\lfloor \frac{j-1}{2} \right\rfloor \right),{\rm \; }j=\overline{1;m_{k} (n)}.\]
 \begin{spacing}{0.3}\end{spacing}
\[b_{j} {\rm \; }\mathop{=}\limits^{\eqref{(6)}} {\rm \; }2\cdot \left(2\cdot \left\lfloor \frac{j}{2} \right\rfloor +j-1-\left\lfloor \frac{j}{2} \right\rfloor \right)=2\cdot \left(j+\left\lfloor \frac{j}{2} \right\rfloor -1\right)=2\cdot \left(j+\frac{j}{2} -\frac{1-(-1)^{j} }{4} -1\right)=\] 
\[=3j-\frac{1-(-1)^{j} }{2} -2=\left(3j+\frac{3-(-1)^{j} }{2} \right)-2+(-1)^{j} -2=\]
 \begin{spacing}{1.1}\end{spacing}
\[=\left(3j+\frac{3-(-1)^{j} }{2} \right)+(-1)^{j} -4=f(j)+(-1)^{j} -4{\rm ,\; }j=\overline{1;m_{k} (n)}. \] 

Therefore, in $k$-row of $T$-matrix, the crossed out elements have form \eqref{(39)}. 

Proposition 5.4 is proved.

PROPOSITION 5.5. For any given natural number $n\ge 5$ the space complexity $C(n)$ of algorithm №1 is
\begin{spacing}{0.5}\end{spacing}
\[\left\lfloor \frac{n+2}{3} \right\rfloor -\left\lfloor \frac{n\% 6}{4} \right\rfloor +\left\lfloor \frac{n\% 6}{5} \right\rfloor .\] 
\begin{spacing}{1.3}\end{spacing}
PROOF. As in the past, let $pr$ is a list comprising numbers of the form $6h\pm 1$ less than or equal to given number $n\in {\rm \mathbb{N}}$: $n\ge 5$. The numbering of the pairs $(6k-1;{\rm \; }6k+1)$ starts with $k=1$, the null element of the list $pr$ is 0. The proof is by induction on number $k\in {\rm \mathbb{N}}$ of pair $(6k-1;{\rm \; }6k+1)$ for each $l\in \{ -1;0;...;4\}$ of $n=6k+l$.

\textbf{Base case. }$k=1$. Therefore, we consider the first pair $(5;7)$. Then, we show that for each $l\in \{ -1;0;...;4\} $ of $n=6+l$ Proposition 5.5 is true. 

$1){\rm \; }l=-1 \Rightarrow C(n)=C(5)=\left\lfloor \frac{5+2}{3} \right\rfloor -\left\lfloor \frac{5\%6}{4} \right\rfloor +\left\lfloor \frac{5\%6}{5} \right\rfloor =\left\lfloor \frac{7}{3} \right\rfloor -\left\lfloor \frac{5}{4} \right\rfloor +\left\lfloor \frac{5}{5} \right\rfloor =2; pr=(0;5).$

$2){\rm \; }l=0 \Rightarrow C(n)=C(6)=\left\lfloor \frac{6+2}{3} \right\rfloor -\left\lfloor \frac{6\%6}{4} \right\rfloor +\left\lfloor \frac{6\%6}{5} \right\rfloor =\left\lfloor \frac{8}{3} \right\rfloor -\left\lfloor \frac{0}{4} \right\rfloor +\left\lfloor \frac{0}{5} \right\rfloor =2; pr=(0;5)).$

$3){\rm \; }l=1 \Rightarrow C(n)=C(7)=\left\lfloor \frac{7+2}{3} \right\rfloor -\left\lfloor \frac{7\%6}{4} \right\rfloor +\left\lfloor \frac{7\%6}{5} \right\rfloor =\left\lfloor \frac{9}{3} \right\rfloor -\left\lfloor \frac{1}{4} \right\rfloor +\left\lfloor \frac{1}{5} \right\rfloor =3; pr=(0;5;7)).$

$4){\rm \; }l=2 \Rightarrow C(n)=C(8)=\left\lfloor \frac{8+2}{3} \right\rfloor -\left\lfloor \frac{8\%6}{4} \right\rfloor +\left\lfloor \frac{8\%6}{5} \right\rfloor =\left\lfloor \frac{10}{3} \right\rfloor -\left\lfloor \frac{2}{4} \right\rfloor +\left\lfloor \frac{2}{5} \right\rfloor =3; pr=(0;5;7)).$

$5){\rm \; }l=3 \Rightarrow C(n)=C(9)=\left\lfloor \frac{9+2}{3} \right\rfloor -\left\lfloor \frac{9\%6}{4} \right\rfloor +\left\lfloor \frac{9\%6}{5} \right\rfloor =\left\lfloor \frac{11}{3} \right\rfloor -\left\lfloor \frac{3}{4} \right\rfloor +\left\lfloor \frac{3}{5} \right\rfloor =3; pr=(0;5;7)).$ 

$6){\rm \; }l=4 \Rightarrow C(n)=C(10)=\left\lfloor \frac{10+2}{3} \right\rfloor -\left\lfloor \frac{10\%6}{4} \right\rfloor +\left\lfloor \frac{10\%6}{5} \right\rfloor =\left\lfloor \frac{12}{3} \right\rfloor -\left\lfloor \frac{4}{4} \right\rfloor +\left\lfloor \frac{4}{5} \right\rfloor =3; pr=(0;5;7).$

\textbf{Inductive hypothesis. }Suppose that Proposition 5.5 is true for any number $k$ of pair
\begin{spacing}{0.5}\end{spacing}
\[(6k-1;6k+1) \text{ at each } l\in \{ -1;0;...;4\} \text{ of } n=6k+l.\]  

\textbf{Inductive step. }Show that Proposition 5.5 is true for $k+1$ at each $l\in \{ -1;0;...;4\} $ of 

\noindent $n=6\cdot (k+1)+l$.
\[\left\lfloor \frac{6\cdot (k+1)+l+2}{3} \right\rfloor -\left\lfloor \frac{(6\cdot (k+1)+l)\% 6}{4} \right\rfloor +\left\lfloor \frac{(6\cdot (k+1)+l)\% 6}{5} \right\rfloor =\] 
\[=\left\lfloor 2+\frac{(6k+l)+2}{3} \right\rfloor -\left\lfloor \frac{(6k+l)\% 6}{4} \right\rfloor +\left\lfloor \frac{(6k+l)\% 6}{5} \right\rfloor =\]
\begin{spacing}{1.2}\end{spacing}
\[=2+\left\lfloor \frac{(6k+l)+2}{3} \right\rfloor -\left\lfloor \frac{(6k+l)\% 6}{4} \right\rfloor +\left\lfloor \frac{(6k+l)\% 6}{5} \right\rfloor .\] 

Then, by the induction hypothesis, we get:
\[\left\lfloor \frac{6\cdot (k+1)+l+2}{3} \right\rfloor -\left\lfloor \frac{(6\cdot (k+1)+l)\% 6}{4} \right\rfloor +\left\lfloor \frac{(6\cdot (k+1)+l)\% 6}{5} \right\rfloor =2+C(k)=C(k+1).\] 

Really,  $l=-1:{\rm \;}\underbrace{0;5;7;...;6k-1;}_{C(k)} {\rm \; }\mathop{}\limits_{+} {\rm \; }\underbrace{6k+1;6\cdot (k+1)-1}_{{\rm 2}} \mathop{}\limits_{6\cdot (k+1)-1} $.
\[l=0:{\rm \; }\underbrace{0;5;7;...;6k-1;}_{C(k)} {\rm \; }\mathop{}\limits_{+} {\rm \; }\underbrace{6k+1;6\cdot (k+1)-1}_{{\rm 2}} \mathop{}\limits_{6\cdot (k+1)} .\] 
\[l=1:{\rm \; }\underbrace{0;5;7;...;6k-1;6k+1}_{C(k)} {\rm \; }\mathop{}\limits_{+} {\rm \; }\underbrace{6\cdot (k+1)-1;6\cdot (k+1)+1}_{{\rm 2}}\mathop{}\limits_{6\cdot (k+1)+1} .\] 
\[l\in \{ 2;3;4\}:{\rm \;}\underbrace{0;5;7;...;6k-1;6k+1}_{C(k)} {\rm \; }\mathop{}\limits_{+} {\rm \; }\underbrace{6\cdot (k+1)-1;6\cdot (k+1)+1}_{{\rm 2}} \mathop{}\limits_{6\cdot (k+1)+l} .\] 

Thus, there is an equation
\begin{spacing}{0.8}\end{spacing}\begin{equation} \label{(42)} 
C(n)=\left\lfloor \frac{n+2}{3} \right\rfloor -\left\lfloor \frac{n\% 6}{4} \right\rfloor +\left\lfloor \frac{n\% 6}{5} \right\rfloor \text{ for all } n\in {\rm \mathbb{N}}: n\ge 5.
\end{equation} 
\begin{spacing}{1.2}\end{spacing}

Proposition 5.5 is proved.

COROLLARY 5.6. $C(n)=O(n)$. 

PROOF. Note that for any natural number $n\ge 5$:

\begin{equation} \label{(43)} 
-\left\lfloor \frac{n\% 6}{4} \right\rfloor +\left\lfloor \frac{n\% 6}{5} \right\rfloor \in \{ -1;0\} .  
\end{equation} 
\[C(n)\mathop{=}\limits^{\eqref{(42)}} \left\lfloor \frac{n+2}{3} \right\rfloor -\left\lfloor \frac{n\% 6}{4} \right\rfloor +\left\lfloor \frac{n\% 6}{5} \right\rfloor \mathop{\le }\limits^{\eqref{(43)}} \left\lfloor \frac{n+2}{3} \right\rfloor \mathop{\le }\limits^{\eqref{(22)}} \frac{n+2}{3}=\]
\begin{spacing}{1.1}\end{spacing}
\[=n-\frac{2n}{3} +\frac{2}{3} =n-\frac{2}{3} \cdot (n-1)<n,{\rm \; }n\ge 5{\rm \; \; }\Rightarrow  {\rm \; \; }C(n)<1\cdot n,{\rm \; }n\ge 5.\] 

It follows that $C(n)=O(n)$. Corollary 5.6 is proved.

\begin{center}
\textbf{\large 6. Conclusion}
\end{center}
\begin{spacing}{1.0}\end{spacing}

Idea of creating $T$-matrix is proposed, basic definitions for her elements introduced. Primary properties of $T$-matrix are explored. Propositions about the form, existence and uniqueness, location of elements in $T$-matrix are proved. Main conclusions, which reflect relationship between $T$-matrix elements and prime numbers, composite numbers $6h\pm 1$, are done. Countability theorems for different sets of $T$-matrix elements, theorems about the «transition down» of $T$-matrix defining element and about the form of singular submatrices that have finite order are proved. Relationship between leading element-counting function $\pi _{{\rm M} _{T} } (x)$ and prime-counting function $\pi (x)$ is established. The upper and lower estimates for $\pi _{{\rm M} _{T} } (x)$ and for leading element $p^{2} (k)$ are got. On the basis of $T$-matrix, prime number generation algorithm for a given natural number $n\ge 5$ is developed. Formulas of composite numbers output from this algorithm at a stages of deletion are found. Asymptotic time complexity and asymptotic space complexity of the algorithm are defined. An exact formula for space complexity of the algorithm is found.

\begin{center}
\textbf{\large References}
\end{center}

\begin{enumerate}

\item Семенов И. Л. Антье и мантисса. Сборник задач с решениями / Под ред. Е. В. Хорошиловой. М.: ИПМ им. М. В. Келдыша, 2015. -- 432 с.

\item  Федоров Ф.М. О рангах и декрементах миноров, определителей и матриц бесконечной системы, 2015.

\item  А.С. Бортаковский, А.В. Пантелеев. Практикум по линейной алгебре и аналитической геометрии: Учеб. пособие -- M.: Высш. шк., 2007. -- 352 с.

\item  Бухштаб А.А. Теория чисел, М., Просвещение, 1966.

\item  Верещагин Н.К., Шень А. Лекции по математической логике и теории алгоритмов. Часть 1. Начала теории множеств. -- 4-е изд., доп. -- М.: МЦНМО, 2012. -- 112 с.

\item  В. И. Зенкин. Распределение простых чисел. Элементарные методы, 2012. -- 112 c.

\item  Ингам А.Э. Распределение простых чисел. М.: Едиториал УРСС, 2005.

\item  Н. Вирт. Алгоритмы и структуры данных. М.: Мир, 1989, 360 стр. 

\item  Чебышев П.Л. О простых числах // Полное собрание сочинений. Т.1. М. -- Л.: АН СССР, 1944. с. 191-207.

\item  J.B. Rosser and L. Schoenfeld. Approximate formulas for some functions of prime numbers, Illinois. J. Math. 6 (1962), 64-94.

\item  J. Sorenson. An introduction to prime number sieves. Technical Report 909, University of Wisconsin, Computer Sciences Department, January 1990.
\end{enumerate}

\end{document}